\documentclass[12pt]{article}
\usepackage{amsmath, amssymb}
\newtheorem{lemma}{Lemma}
\newtheorem{corollary}{Corollary}
\newtheorem{proposition}{Proposition}
\newtheorem{definition}{Definition}
\newtheorem{conjecture}{Conjecture}

\newcommand{\adl}{{\mbox{\upshape{ad}}_l}}
\newcommand{\adr}{{\mbox{\upshape{ad}}_r}}
\newcommand{\adrtil}{\widetilde{{\mbox{\upshape{ad}}_r}}}
\newcommand{\Bc}{{\check{B}}}
\newcommand{\Btil}{{\tilde{B}}}
\newcommand{\C}{{\mathbb C}}

\newcommand{\N}{{\mathbb N}}

\newcommand{\cC}{{\mathcal C}}
\newcommand{\Cc}{{\mathcal C}}

\newcommand{\cI}{{\mathcal I}}
\newcommand{\cIm}{{\mathcal I}_\mu}

\newcommand{\cM}{{\mathcal M}}

\newcommand{\cop}{{\mbox{\scriptsize{cop}}}}

\newcommand{\End}{\mbox{End}}

\newcommand{\gfrak}{{\mathfrak g}}

\newcommand{\hfrak}{{\mathfrak h}}

\newcommand{\Ibar}{\overline{I}}
\newcommand{\id}{{\mbox{id}}}

\newcommand{\Jtil}{{\tilde{J}}}
\newcommand{\kb}{\kappa_\Bc}
\newcommand{\kow}{{\varDelta}}

\newcommand{\ltil}{\tilde{l}}

\newcommand{\mutil}{\tilde{\mu}}

\newcommand{\op}{{\mbox{\scriptsize{op}}}}
\newcommand{\ot}{\otimes}
\newcommand{\Ppi}{P(\pi)}
\newcommand{\Ppiplus}{P^+(\pi)}
\newcommand{\pt}{{\pi_\Theta}}
\newcommand{\pti}{{\pi_{\Theta,i}}}
\newcommand{\Ptheta}{P_{Z(\Bc)}}
\newcommand{\qfield}{\cC}
\newcommand{\rf}{{\bf{r}}}
\newcommand{\rfbar}{{\bf{\bar{r}}}}
\newcommand{\rk}{\mbox{rank}}

\newcommand{\rtil}{\tilde{R}}

\newcommand{\ract}{\triangleleft}
\newcommand{\rqg}{R_q[G]}
\newcommand{\rqbp}{R_q[B^+]}
\newcommand{\rqbm}{R_q[B^-]}

\newcommand{\thetatil}{\tilde{\theta}}

\newcommand{\Ttcp}{{T_\Theta '}}

\newcommand{\uqbpc}{{\check{U}_q(\mathfrak{b}^+)}}
\newcommand{\uqbmc}{{\check{U}_q(\mathfrak{b}^-)}}
\newcommand{\uqnp}{{U_q(\mathfrak{n}^+)}}
\newcommand{\uqnm}{{U_q(\mathfrak{n}^-)}}

\newcommand{\uqg}{{U_q(\mathfrak{g})}}
\newcommand{\uqgc}{{\check{U}_q(\mathfrak{g})}}
\newcommand{\Tc}{{\check{T}}}
\newcommand{\Uc}{{\check{U}}}
\newcommand{\vep}{\varepsilon}
\newcommand{\wght}{\mathrm{wt}}
\newcommand{\wurz}{\pi}
\newcommand{\Xtil}{{\tilde{X}}}
\newcommand{\Ytil}{{\tilde{Y}}}
\newcommand{\Z}{{\mathbb Z}}

\title{Quantum Symmetric Pairs and the Reflection Equation}

\author{Stefan Kolb\footnote{Supported by the German Research Foundation (DFG)}}

\date{August 21, 2006}

\begin{document}

\maketitle

\begin{abstract}
  It is shown that central elements in G.~Letzter's quantum group analogs of symmetric pairs 
  lead to solutions of the reflection equation. This clarifies the relation between 
  Letzter's approach to quantum symmetric pairs and the approach taken by M.~Noumi,
  T.~Sugitani, and M.~Dijkhuizen. 
  
  We develop general tools to show that a Noumi-Sugitani-Dijkhuizen type construction 
  of quantum symmetric pairs can be performed preserving spherical 
  representations from the classical situation. 
  These tools apply to the symmetric pair $FII$ and to all symmetric pairs 
  which correspond to an automorphism of the underlying Dynkin diagram. 
  Hence Noumi-Sugitani-Dijkhuizen type constructions with desirable properties are possible 
  for various symmetric pairs for exceptional Lie algebras. 
\end{abstract}

\section{Introduction}
The reflection equation with spectral parameter first appeared in I.~V.~Cherednik's work
on factorizing scattering on the half line \cite{a-Cher84}. The present paper is devoted
to the reflection equation without spectral parameter \cite{a-KSS93} and its relation
to quantum group analogs of symmetric pairs.

Let $\theta$ be an involutive automorphism of a complex semisimple Lie algebra $\gfrak$ 
and let $\gfrak^\theta$ denote the fixed Lie subalgebra. The pair $(\gfrak,\gfrak^\theta)$
is a classical (infinitesimal) symmetric pair. 

In a series of papers \cite{a-Noumi96}, \cite{a-NS95}, \cite{a-Dijk96}, \cite{a-DN98},
\cite{a-DS99}, and references therein, analogs of $U(\gfrak^\theta)$ inside the simply 
connected quantum universal
enveloping algebra $\uqgc$ were constructed for all classical symmetric pairs. The 
constructions were done case by case and depend on solutions of the 
reflection equation. The solutions of the reflection equation were
given explicitly. 

Taking a unifying approach, G.~Letzter classified all maximal left coideal 
subalgebras $\Bc\subseteq\uqgc$ which specialize to $U(\gfrak^\theta)$ as $q$ goes to 1
\cite{MSRI-Letzter}. These coideal subalgebras can be given explicitly in terms of 
generators and relations \cite{a-Letzter03}. However, the reflection equation does not 
appear in this construction. 
In \cite[Section 6]{a-Letzter99a} it was shown that the coideal subalgebras constructed via
solutions of the reflection equation are contained in Letzter's coideal subalgebras. The 
proof depends on the abstract characterization of quantum symmetric pairs. 

In the present paper we consider the converse problem: Let $\Bc\subseteq \uqgc$ be one of 
Letzter's coideal subalgebras. 
\begin{align*}
  \mbox{ \begin{minipage}{13.5cm} 1) Is there a solution of the reflection equation 
  canonically associated to $\Bc$ which can be used for a 
  Noumi-Sugitani-Dijkhuizen type construction of a subcoideal $L\subseteq \Bc$ ?
  \end{minipage}}
\end{align*}
In order to replace $\Bc$ by $L$ when looking at invariants, we moreover ask the following 
question:
\begin{align*}
  \mbox{ \begin{minipage}{13.5cm} 2) Let $V$ be a type one representation of $\uqgc$.
  Does the space of $L$-invariant elements in $V$ coincide with the subspace of $\Bc$-invariant
  elements ?
  \end{minipage}}
\end{align*}
In this paper question 1) is answered in the affirmative. It is shown that suitable elements 
in the center $Z(\Bc)$ lead to solutions of the reflection equation. The main ingredient in 
the argument is an isomorphism between the $q$-deformed coordinate ring $\rqg$ and the 
locally finite part of $\uqgc$. This isomorphism is given by $l$-functionals and is the
inverse of Caldero's isomorphism \cite{a-Caldero93} given by the Rosso form.
The center $Z(\Bc)$ is calculated in \cite{a-KLp} where it is shown that the desired central
elements exist. The knowledge of the center together with a suitable notion of 
minimal weight leads to solutions of the reflection equation also for exceptional 
Lie algebras.

The main motivation for question 1) is to gain a better understanding
of the ring of invariants $\rqg^\Bc=\{x\in \rqg\,|\, bx=\vep(b)x\,\forall b\in \Bc\}$. 
For the quantum symmetric pairs of classical type the solution of the reflection equation
was used to obtain sets of generators of $\rqg^\Bc$ \cite{a-Noumi96}, \cite{a-DN98}. It is planned to investigate a similar
construction for general quantum symmetric pairs. 

An affirmative answer to Question 2) on the other hand would show that all results about
zonal spherical functions and Macdonald polynomials \cite{a-Letzter03}, \cite{a-Letzter04} 
also hold for the coideal $L$. 
In the present paper this is proved for the symmetric pair $FII$. 
We also formulate a generalization of Question 2) where $L$ is replaced by the left
coideal generated by any nontrivial element of the basis of the center $Z(\Bc)$ determined
in \cite{a-KLp}. This generalization is proved for all symmetric pairs which correspond to 
automorphisms of the underlying Dynkin diagram. 

The paper is organized as follows. In Section \ref{QG} we fix notations and recall
general facts about the universal $r$-form for $\rqg$. These facts are used to 
show that $l$-functionals give an isomorphism between $\rqg$ and the locally finite part of 
$\uqg$. 

In Section \ref{Coideal-Refl} it is shown how suitable central elements of a coideal
subalgebra of $\uqgc$ lead to solutions of the reflection equation. This result
is applied to Letzter's family of quantum symmetric pairs in \ref{invariants}.
To this end Letzter's construction of quantum symmetric pair coideal subalgebras and
properties of their center are recalled in Sections \ref{qsp} and \ref{cqsp}, respectively.
Moreover, in \ref{minweights} a natural notion of minimal weights is introduced.
A precise formulation of Question 2) above is given in Section \ref{invariants} as 
Conjecture~\ref{conj1} and a formulation of its generalization is given as 
Conjecture~ \ref{conj2}. 

In this paper, following \cite{MSRI-Letzter} we consider left coideal subalgebras and
the right adjoint action. In \ref{leftright} the reflection equation obtained here is 
translated into the setting of right coideal subalgebras reproducing  precisely the 
reflection equation \cite[(5.5)]{a-Dijk96}, \cite[(6.4)]{a-DS99}.

In Section \ref{towards} we develop general results about the left coideal generated
by a nontrivial element of the basis of $Z(\Bc)$. These results allow us to prove
Conjecture \ref{conj1} for the symmetric pair of type $FII$. Section \ref{pt-empty}
is devoted to symmetric pairs which correspond to an automorphism of the 
underlying Dynkin diagram. In this case weight considerations and a suitable projection
map allow us to recover all generators of $\Bc$ inside the coideal generated by an 
element of the basis of $Z(\Bc)$. This proves Conjecture \ref{conj2} in this special case.

The author is very grateful to Gail Letzter for encouragement and many instructive 
discussions about quantum symmetric pairs. He also wishes to thank Jasper Stokman for 
his interest and helpful comments.

\section{Quantum groups}\label{QG}
We write $\C$, $\Z$, and $\N_0$ to denote the complex numbers, the integers, and the nonnegative
integers, respectively.
  \subsection{Notations}\label{notation1}
    Let $\gfrak$ be a finite dimensional complex semisimple 
    Lie algebra of rank $n$ and $\hfrak\subseteq\gfrak$ a fixed Cartan subalgebra.
    Let $\Delta\subseteq \hfrak^\ast$ denote the root system
    associated with $(\gfrak,\hfrak)$.
    Choose an ordered basis $\wurz=\{\alpha_1,\dots,\alpha_n\}$ of simple 
    roots for $\Delta$ and let $\Delta^+$ (resp.~$\Delta^-$) be the set of positive 
    (resp.~negative) roots with respect to $\wurz$.
    Identify $\hfrak$ with its dual via the Killing form. The induced
    nondegenerate symmetric bilinear form on $\hfrak^*$
    is denoted by $(\cdot,\cdot)$. 
    We write $Q(\pi)$ for the root lattice and $P(\pi)$ for the weight lattice
    associated to the root system $\Delta$. Let
    $\omega_i\in\hfrak^\ast$, $i=1,\dots,n$ be the fundamental weights with respect to $\pi$
    and $P^+(\wurz)$ denote the set of dominant weights, i.~e.~the 
    $\N_0$-span of $\{\omega_i\,|\,i=1,\dots,n\}$.
    Let $\leq$ denote the standard partial ordering on 
    $\hfrak^*$.  In particular, $\mu\leq\gamma$ if and only  if $\gamma-\mu\in \N_0\pi$.
    Finally, let $W$ denote the Weyl group associated to the root system $\Delta$ and let
    $w_0\in W$ be the longest element in $W$ with respect to $\pi$. 
  
  \subsection{$\uqgc$ and $\rqg$}\label{uqgrqg}
    Let $\qfield=\C(q^{1/N})$ denote the field of rational functions in one
    variable $q^{1/N}$ where $N$ is sufficiently large such that 
    $(\lambda,\mu)\in \frac{1}{N} \Z$ for all $\lambda,\mu\in \Ppi$.
    We consider here the simply connected quantum universal enveloping algebra 
    $\uqgc$ as the $\qfield$-algebra generated by
    elements $\{x_i, y_i, \tau(\lambda)\,|\, i=1,\dots,n,\,\lambda\in \Ppi\}$ and 
    relations as given for instance in \cite[Section 3.2.9]{b-Joseph}. 
    As usual we will write $t_i=\tau(\alpha_i)$. For $\alpha\in Q(\wurz)$ and any subset 
    $M\subseteq \uqgc$ let
    \begin{align}\label{wspace}
      M_\alpha:=\{u\in M\,|\, \tau(\lambda)u\tau(-\lambda)=q^{(\lambda,\alpha)}u 
      \mbox{ for all $\lambda \in P(\pi)$}\}
    \end{align}  
    denote the set of elements of weight $\alpha$ in $M$.
    The algebra $\uqgc$ has a Hopf algebra structure with counit $\vep$, coproduct 
    $\kow$, and antipode $\sigma$ as in \cite[Section 3.2.9]{b-Joseph}.
    In particular the coproduct satisfies
    \begin{align}\label{kopr}
      \kow x_i=&t_i\ot x_i + x_i\ot 1,\qquad \kow y_i=1\ot y_i + y_i\ot t_i^{-1},\qquad
      \kow t_i=t_i\ot t_i.
    \end{align}
    Recall that $\sigma^2$ is a Hopf algebra automorphism of $\uqg$ and that
      \begin{align}\label{squareofAnti}
        \sigma^2(a)=\tau(-2\rho)a\tau(2\rho)
      \end{align}
    where $\rho=\frac{1}{2} \sum_{\alpha\in \Delta^+}\alpha$ is half the sum 
    of all positive roots of $\gfrak$. 
    
    As usual let $\uqnp$, $\uqnm$, $G^+$, and $G^-$ denote the subalgebras of $\uqgc$ 
    generated by the sets $\{x_i\,|\, i=1,\dots,n\}$, $\{y_i\,|\, i=1,\dots,n\}$, 
    $\{x_it_i^{-1}\,|\, i=1,\dots,n\}$, and $\{y_it_i\,|\, i=1,\dots,n\}$, respectively.
    Let moreover $\uqbpc$ and $\uqbmc$ denote the subalgebras of $\uqgc$ generated by the 
    sets $\{x_i,\tau(\lambda)\,|\, i=1,\dots,n,\lambda \in \Ppi\}$ and 
    $\{y_i,\tau(\lambda)\,|\, i=1,\dots,n,\lambda \in \Ppi\}$, respectively. We write
    $\Tc$ for the multiplicative group generated by $\tau(\lambda)$, $\lambda\in P(\pi)$,
    and $\Uc^0$ for the subalgebra of $\uqgc$ generated by $\Tc$.
    The algebra $\uqgc$ admits a triangular decomposition. More precisely, the multiplication
    map induces an isomorphism
    \begin{align*}
      \uqnm\ot \Uc^0 \ot \uqnp \rightarrow \uqgc.
    \end{align*}
    The same holds if one or both of $\uqnm$ and $\uqnp$ are replaced by $G^-$ and $G^+$,
    respectively. 
    
    Note that $\uqnp$, $\uqnm$, $G^+$, and $G^-$ are direct sums of their weight spaces.
    Using this decomposition into weight spaces and the triangular decomposition one can 
    define projections
    \begin{align}\label{pialphbet}
      \pi_{\alpha,\beta}:\uqgc\cong \uqnm\ot \Uc^0\ot G^+ \rightarrow
                    \uqnm_\alpha\ot \Uc^0\ot G^+_\beta
    \end{align}
    for all $\alpha,\beta\in Q(\pi)$. Moreover, define
    \begin{align}\label{pialph}
      \pi_\alpha:\uqgc\rightarrow \uqgc_\alpha
    \end{align}
    to be the projection onto the weight space of weight $\alpha\in Q(\pi)$.
    
    The following notation is used in Sections \ref{qsp} and \ref{pt-empty}.
    For any multiindex $I=(i_1,\dots,i_m)$, $1\le i_j\le n$, define
    $y_I=y_{i_1}\cdots y_{i_m}$ and $x_I=x_{i_1}\cdots x_{i_m}$.
    Set ${\rm wt}(I)=\alpha_{i_1}+\dots+\alpha_{i_m}$. 
    Let $\cI$ be a set of multiindices such that $\{y_I\,|\,I\in \cI\}$
    is a basis of $\uqnm$.
    
    For $\mu\in \Ppiplus$ let $V(\mu)$ denote the uniquely determined
    finite dimensional simple left $\uqgc$-module of highest weight 
    $\mu$. More explicitly, there exists a highest weight vector
    $v_\mu\in V(\mu)\setminus\{0\}$ such that
    \begin{align}\label{hw}
      x_i v_\mu=0,\qquad \tau(\lambda) v_\mu=q^{(\mu,\lambda)}v_\mu \qquad
       \forall\, i=1,\dots,n\,\mbox{ and }\, \lambda \in \Ppi.
    \end{align}
    For any weight vector $v\in V(\mu)$ we will write $\wght(v)$ to denote its weight.
    We consider the dual space $V(\mu)^*$ always with its natural right $\uqgc$-module 
    structure defined by $v^\ast u(v)=v^\ast (uv)$ for all $v^\ast \in V(\mu)^*$, 
    $v\in V(\mu)$, 
    and $u\in \uqgc$. We call an element $v^\ast$ of the right $\uqgc$-module $V(\mu)^\ast$
    a weight vector of weight $\nu\in P(\pi)$ if 
    $v^\ast \tau(\lambda)=q^{(\nu,\lambda)}v^\ast$ for
    all $\lambda\in P(\pi)$. Moreover, we write $V(\mu)^\ast_\nu$ to denote the space spanned
    by all weight vectors of weight $\nu$ in $V(\mu)^\ast$.
    
    The $q$-deformed coordinate ring $\rqg$ is defined  \cite[9.1.1]{b-Joseph} to be the 
    subspace of the linear dual $\uqg^\ast$ spanned by the matrix coefficients 
    of the finite dimensional irreducible representations $V(\mu)$, $\mu\in \Ppiplus$.
    For $v\in V(\mu)$, $v^\ast \in V(\mu)^\ast$ the matrix coefficient
    $c^\mu_{v^\ast,v}\in \uqg^\ast$ is defined by
        \begin{align*}
          c^\mu_{v^\ast,v}(X)=v^\ast(Xv).
        \end{align*}
    The linear span of matrix coefficients of $V(\mu)$ will be denoted by $C^{V(\mu)}$.
    Counit, antipode, and coproduct of $\rqg$ will also be denoted by $\vep$, $\sigma$,
    and $\kow$, respectively. Recall that $\rqg$ is a $\uqgc$-bimodule, with bimodule 
    structure defined by $(uav)(x)=a(vxu)$ for all $a\in \rqg$, $u,v,x\in\uqgc$.
    
    Let $H$ be a coalgebra with coproduct $\kow$. We will frequently make use of 
    Sweedler notation in the form $\kow h=h_{(1)}\ot h_{(2)}$. Moreover, we write $H^\cop$
    to denote the coalgebra with the opposite coproduct $\kow^\cop h=h_{(2)}\ot h_{(1)}$.
    Similarly, for any algebra $A$ we write $A^\op$ to denote the algebra with the opposite
    product.
    
  \subsection{The universal $r$-form for $\rqg$}
    The existence of the Rosso form for $\uqg$ is reflected in the fact that $\rqg$ is
    a coquasitriangular Hopf algebra. Recall (e.g.~\cite[chapter 10]{b-KS})  that a 
    coquasitriangular bialgebra over $\qfield$ is a bialgebra $A$ over $\qfield$ equipped 
    with a convolution invertible skew-pairing $\rf: A\ot A\rightarrow \qfield$ such that
    \begin{align*}
      m_{A^\op}=\rf\ast m_A\ast \rfbar. 
    \end{align*}
    Here $m_A:A\ot A\rightarrow A$ denotes the multiplication map of 
    the algebra $A$, the symbol $\ast$ denotes the convolution product
    and $\rfbar$ is the convolution inverse of $\rf$.
    If $(A,\rf)$ is a coquasitriangular bialgebra then $\rf$ is called a universal $r$-form for
    $A$. We quickly review the construction of the universal $r$-form for $\rqg$ along the
    lines of \cite{a-Gaitsgory95}, \cite{b-Joseph} in our notational conventions.
       
    Let $\rqbp$ and $\rqbm$ denote the restriction of 
    $\rqg\subseteq\uqgc^\ast$ to $\uqbpc$ and $\uqbmc$, respectively. 
    The Hopf algebra structure of $\rqg$ induces Hopf algebra
    structures on $\rqbp$ and $\rqbm$. Let 
    $\Psi^+:\rqg\rightarrow \rqbp$ and $\Psi^-:\rqg\rightarrow \rqbm$
    denote the canonical projections.
    
    By \cite[Section 2]{a-Tanisaki92} there exists a
    unique  skew-pairing $\varphi:\uqbmc\times\uqbpc\rightarrow \qfield$
    such that
    \begin{align*}
      \varphi(y_i,x_j)&=-\delta_{ij}/(q^{(\alpha_i,\alpha_i)/2}-q^{-(\alpha_i,\alpha_i)/2})\\
      \varphi(\tau(\lambda),\tau(\mu))&=q^{-(\lambda,\mu)}.
    \end{align*}    
    By \cite[9.2.11, 9.2.12]{b-Joseph} the pairing $\varphi$ can be used to define  
    Hopf algebra isomorphisms
    \begin{align*}
          \Phi^-:\uqbpc^{\op}\rightarrow \rqbm, 
                       \qquad x\mapsto \varphi(\cdot,x)\\
          \Phi^+:\uqbmc^{\cop}\rightarrow \rqbp, 
                       \qquad y\mapsto \varphi(y,\cdot).            
    \end{align*}
    The universal $r$-form for $\rqg$ is defined to be the composition of maps
       \begin{align*}
          \rf:\rqg\ot\rqg \stackrel{\Psi^+\ot\Psi^-}{\longrightarrow}
          &\rqbp\ot\rqbm\\
          \stackrel{(\Phi^+)^{-1}\ot (\Phi^-)^{-1}}{\longrightarrow}
          &\uqbmc\ot\uqbpc\stackrel{\varphi}{\longrightarrow}\qfield.
        \end{align*}
    By \cite[2.2.3]{a-Gaitsgory95}, \cite[9.4.7]{b-Joseph} the map $\rf$ is a universal 
    $r$-form for $\rqg$.
    For the purpose of the present paper the following properties of $\rf$ are relevant.
        \begin{proposition}\label{r-trig}
          Let $\lambda,\mu\in \Ppiplus$ and let $\{v_i\}\subseteq V(\lambda)$
          and $\{w_j\}\subseteq V(\mu)$ be weight bases with dual bases
          $\{v_i^\ast\}\subseteq V(\lambda)^\ast$ and 
          $\{w_j^\ast\}\subseteq V(\mu)^\ast$, respectively. Then the following
          properties hold.
          
          i) $\rf(c^\lambda_{v_i^\ast,v_{i'}},c^\mu_{w_j^\ast,w_{j'}})\neq 0
             \,\Longrightarrow \wght(v_i)\ge \wght (v_{i'})\,\mbox{and }
             \wght(w_j)\le \wght (w_{j'})$.
   
          ii) $\rf(c^\lambda_{v_i^\ast,v_{i'}},c^\mu_{w_j^\ast,w_{j'}})\neq 0
             \,\Longrightarrow \wght(v_i)-\wght (v_{i'})= \wght(w_{j'})- \wght (w_j)$.
   
          iii) $\rf(c^\lambda_{v_i^\ast,v_{i}},c^\mu_{w_j^\ast,w_{j}})
               =\varphi(\tau(\wght(v_i)),\tau(\wght(w_j)))=
               q^{-(\wght(v_i),\wght(w_j))}$.
               
          iv) $\rf(\sigma(a),\sigma(b))=\rf(a,b)$.     
        \end{proposition}
        {\bf Proof:} Property $i)$ follows from the fact that $\rf$ factors
        over $\rqbp\ot \rqbm$. Property $ii)$ follows from the properties of $\varphi$
        \cite[Lemma 2.1.3]{a-Tanisaki92}.
        To verify property $iii)$ note that
        \begin{align*}
          \Phi^+(\tau(-\wght(v_i)))&=c^\lambda_{v_i^\ast,v_{i}}\big|_{\rqbp},\\
          \Phi^-(\tau(-\wght(w_j)))&=c^\mu_{w_j^\ast,w_{j}}\big|_{\rqbm}.
        \end{align*}  
        The last statement follows from \cite[Lemma 2.1.2]{a-Tanisaki92}. 
        $\blacksquare$

  \subsection{The locally finite part}\label{locfinpart}
    Any Hopf algebra $H$ with antipode $\sigma$ acts on itself from the left and the right 
    by the left and right adjoint action, respectively. Using Sweedler notation these actions
    are given explicitly as follows
    \begin{align*}
      \adl(h)u=h_{(1)}u \sigma(h_{(2)}),\qquad \adr(h)u=\sigma(h_{(1)})u h_{(2)},\qquad
      u,h\in H.
    \end{align*}
    The right locally finite part is defined by
    \begin{align*}
      F_r(H)=\{h\in H\,|\, \dim((\adr H)h)<\infty \}.
    \end{align*}
    Similarly, the left locally finite part $F_l(H)$ is defined using
    the left adjoint action $\adl$ instead of $\adr$. Note that 
    for $u,x\in H$ on has $(\adr u)x=\sigma((\adl(\sigma^{-1}(u))\sigma^{-1}(x))$
    and therefore
    \begin{align}\label{left-right-locfin}
      F_r(H)=\sigma(F_l(H)).
    \end{align}
    It was shown in \cite[Thm.~4.10]{a-JoLet2}, \cite{a-Caldero93} that the locally 
    finite part of $\uqgc$ can be written explicitly as
    \begin{align*}
       F_l(\uqgc)=\bigoplus_{\mu\in\Ppiplus}(\adl \uqgc)\tau(-2\mu).
    \end{align*} 
    By (\ref{left-right-locfin}) this implies
    \begin{align*}
           F_r(\uqgc)=\bigoplus_{\mu\in\Ppiplus}(\adr \uqgc)\tau(2\mu).
    \end{align*} 
    The following observation will be frequently used in Section \ref{towards}.
    For any Hopf algebra $H$ the right adjoint action is compatible with the coproduct in 
    the following sense
    \begin{align}\label{kowadr}
      \kow(\adr b) C=&(\adr b_{(2)})C_{(1)}\ot \sigma(b_{(1)})C_{(2)}b_{(3)}\qquad
      \mbox{for all $b,C\in H$}
    \end{align}  
    In particular in view of (\ref{kopr}) one obtains for the generators $x_i,y_i$ of $\uqgc$
    the relations
    \begin{align}            
      \kow(\adr x_i) C=&(\adr x_i)C_{(1)}\ot t_i^{-1}C_{(2)} + 
                        (\adr t_i)C_{(1)}\ot t_i^{-1}C_{(2)} x_i \nonumber\\
                        &-C_{(1)}\ot t_i^{-1} x_i C_{(2)},\label{adrxC}\\
      \kow(\adr y_i) C=&(\adr y_i)C_{(1)}\ot C_{(2)}t_i^{-1} + C_{(1)}\ot C_{(2)} y_i\nonumber\\
                               &  -(\adr t_i^{-1})C_{(1)}\ot y_i t_i C_{(2)} t_i^{-1}. \label{adryC}              
    \end{align} 
    for all $C\in \uqgc$.
    
    \medskip
    \noindent
    {\bf Remark:} For $H=\uqgc$ it was proved explicitly \cite[5.3]{a-JoLet95},
        \cite[7.1.6]{b-Joseph}, \cite[Thm. 5.1]{MSRI-Letzter} that
        $F_l(H)$ is a left coideal of $H$. However, this result holds for
        any Hopf algebra $H$. It can be proved analogously to \cite[Lemma~5]{a-heko03}.
        
\subsection{$l$-functionals}    
  We now apply the theory of $l$-functionals (e.g.~\cite[10.1.3]{b-KS}) to the
  coquasitriangular Hopf algebra $\rqg$. This will lead to an isomorphism between the
  matrix coefficients $C^{V(\mu)}$ and the corresponding component of the locally finite
  part. 
  
  Using the universal $r$-form one defines $l$-functionals in the dual 
  Hopf algebra $\rqg^\circ$ by
  \begin{align*}
    l^+(a):=\rf(\cdot,a)&,\quad  l^-(a):=\rf(\sigma(a),\cdot),\\
    l(a):=l^-(\sigma^{-1}&(a_{(1)})) l^+(a_{(2)})
          =\rf(a_{(1)},\cdot) \rf(\cdot,a_{(2)}),\\    
    \ltil(a):=l^+(a_{(1)})&l^-(\sigma^{-1}(a_{(2)}))
         =\rf(\cdot,a_{(1)})\rf(a_{(2)},\cdot)
  \end{align*}
  for $a\in \rqg$.
  Note that by definition
  \begin{align*}
    l^+(a)&=(\Phi^-)^{-1}\circ\Psi^-(a)\in \uqbpc,\\
    l^-(a)&=(\Phi^+)^{-1}\circ\Psi^+(\sigma(a))\in \uqbmc
  \end{align*}
  and therefore $l(a), \ltil(a)\in \uqgc$ and
  \begin{align*}
    \rf(a,b)=l^+(b)(a)=\varphi(l^-(\sigma^{-1}(a)),l^+(b)).
  \end{align*}
  Using the fact that the universal $r$-form is a skew-pairing one
  can determine the coproduct of the $l$-functionals. 
  One obtains \cite[10.1.3, Prop.~11]{b-KS}
  \begin{align}
    &\kow l^+(a)=l^+(a_{(1)})\ot l^+(a_{(2)}),\quad
    \kow l^-(a)=l^-(a_{(1)})\ot l^-(a_{(2)}),\nonumber\\
    &\kow l(a)=l(a_{(2)})\ot l^-(\sigma^{-1}(a_{(1)})) l^+(a_{(3)}) 
         =l(a_{(2)})\ot \sigma(l^-(a_{(1)}))  l^+(a_{(3)}),\label{lcoprod}\\
   &\kow \ltil(a)=l^+(a_{(1)}) l^-(\sigma^{-1}(a_{(3)}))\ot \ltil(a_{(2)})
           =l^+(a_{(1)}) \sigma(l^-(a_{(3)}))\ot \ltil(a_{(2)}) .\nonumber 
  \end{align}
  The $l$-functionals $l$ and $\ltil$ are compatible with the right and left adjoint action,
  respectively, in the following way \cite[10.1.3, Prop.~11]{b-KS}
  \begin{align}
     \adr(u)(l(a))&=(a_{(1)}\sigma(a_{(3)}))(u)\,l(a_{(2)})\label{adrcompatible}\\
     \adl(u)(\ltil(a))&=(\sigma(a_{(1)})a_{(3)})(u)\,\ltil(a_{(2)}).\label{adlcompatible}
  \end{align}  
  These formulas imply in particular that for all $a\in \rqg$
  \begin{align*}
    l(a)\in F_r(\uqgc),\qquad \ltil(a)\in F_l(\uqgc).
  \end{align*}  
  The following proposition is a generalization of a result proved for $\mathfrak{sl}(r+1)$ 
  and
  $\mathfrak{sp}(2r)$ in \cite[Lemma 4.7]{a-HeckSchm2}. The proof given here is significantly simpler.
  
  \begin{proposition}\label{l-iso}
    For any $\mu\in \Ppiplus$ the $l$-functionals $l$ and $\ltil$ define
    isomorphisms of vector spaces
    \begin{align}\label{l-map}
      l:C^{V(\mu)}\rightarrow (\adr \uqgc)\tau(-2w_0\mu),\qquad
      \ltil:C^{V(\mu)}\rightarrow (\adl \uqgc)\tau(-2\mu)
    \end{align}
    Moreover, $l$ and $\ltil$ are compatible with the right and left adjoint action,
    respectively, in the following way
    \begin{align}\label{l-compatible}
      l(\sigma(u_{(2)})au_{(1)})=(\adr u)l(a), \quad  
      \ltil(u_{(2)}a\sigma(u_{(1)}))=(\adl u)\ltil(a),
    \end{align}
    for all $u\in \uqgc$, $a\in \rqg$. 
  \end{proposition}
  {\bf Proof:} Formulas (\ref{l-compatible}) follow immediately from
  (\ref{adrcompatible}) and (\ref{adlcompatible}). 
  
  For any dominant integral weight $\mu\in \Ppiplus$ let 
  $v_{w_0\mu}\in V(\mu)$ denote a lowest weight vector and 
  let $v^\ast_{w_0\mu}\in V(\mu)^\ast$ be a lowest weight vector 
  in the dual right $\uqgc$-module such that $v^\ast_{w_0\mu}(v_{w_0\mu})=1$. 
  In $\uqgc$ using Proposition \ref{r-trig} one verifies
  \begin{align}\label{get-tau2lambda}
    \tau(-2w_0\mu)=l(c^\mu_{v^\ast_{w_0\mu},v_{w_0\mu}}).
  \end{align}
  To shorten notation define 
  $c:=c^\mu_{v^\ast_{w_0\mu},v_{w_0\mu}}$ and define a
  right action $\adrtil$ of $\uqgc$ on $\rqg$ by
  \begin{align*}
    (\adrtil u)a=\sigma(u_{(2)})au_{(1)}.
  \end{align*}
  Then by (\ref{l-compatible}) and (\ref{get-tau2lambda}) the map $l$ 
  maps $(\adrtil \uqgc)c$ onto the right $\uqgc$-module
  $(\adr(\uqgc))(\tau(-2w_0\mu))$.
  Moreover, by \cite[Thm.~3.5]{a-JoLet2} one has
  \begin{align*}
    \dim((\adr \uqgc)\tau(-2w_0\mu))=
  \dim(\End(V(\mu)))=\dim (C^{V(\mu)})
  \end{align*}
  and therefore $(\adrtil \uqgc)c=C^{V(\mu)}$ and 
  $l|_{C^{V(\mu)}}$ is also injective. This proves the first isomorphism in 
  (\ref{l-map}). The second isomorphism is obtained analogously. 
  $\blacksquare$
  
  \vspace{.3cm}
  
  \noindent{\bf Remark:} P.~Caldero showed that the Rosso form induces an isomorphism
  $(\adl \uqgc)\tau(-2\mu)\rightarrow C^{V(\mu)}$ \cite{a-Caldero93}.
  Using the properties of $\ltil$ and of the Rosso form one can verify that $\ltil$
  is the inverse of this isomorphism.
  
  \section{Reflection equation for quantum symmetric pairs}\label{Reflection}
  
  \subsection{Coideal subalgebras and a reflection equation}\label{Coideal-Refl}
  Let $B\subseteq \uqgc$ be a left coideal subalgebra with center $Z(B)$. In the previous 
  section we saw that $(\adr \uqgc)\tau(-2w_0\mu)$ is spanned by $l$-functionals. This
  can now be used to show that any element in $Z(B)\cap (\adr \uqgc)\tau(-2w_0\mu)$ gives 
  rise to a solution of the reflection equation. 
  
  For any $\mu \in \Ppiplus$ let $N=\dim V(\mu)$ and let 
  $\{v_1,\dots,v_N\}$ be a basis
  of $V(\mu)$ with dual basis $\{v_1^\ast,\dots,v^\ast_N\}$.
  To shorten notation define $c^i_j:=c^\mu_{v_i^\ast,v_j}$ for 
  $i,j=1,\dots,N$ and
  \begin{align*}
    R^{ij}_{kl}:=\rf(c^i_k,c^j_l).
  \end{align*}
  By Proposition \ref{l-iso} any element $C\in (\adr \uqgc)\tau(-2w_0\mu)$ can be written 
  in the form
  \begin{align}\label{J-def}
    C=\sum_{i,j}J^j_i l(c^i_j)
  \end{align}
  for some coefficients $J^j_i\in \qfield$. 
  \begin{proposition}\label{reflprop}
    Given $\mu\in \Ppiplus$, a left coideal subalgebra 
    $B\subseteq \uqgc$ and an element
    \begin{align*}
      C\in Z(B)\cap(\adr \uqgc)\tau(-2w_0\mu).
    \end{align*}
    Then the matrix $J=(J^j_i)_{i,j=1,\dots,N}$ 
    defined by (\ref{J-def}) satisfies the reflection equation
    \begin{align}
      J_1\rtil J_2 R = R_{21} J_2 \rtil_{21} J_1 \label{refl}
    \end{align} 
    where $\rtil^{ij}_{kl}=\rf(c^i_k,\sigma(c^j_l))$,
    the lower index $21$ denotes conjugation with the twist 
    $P:v\ot w\mapsto w\ot v$, and ${J_1}^{ij}_{kl}=\delta_{jl}J^i_k$,
    ${J_2}^{ij}_{kl}=\delta_{ik}J^j_l$. 
  \end{proposition}
  {\bf Proof:}
   As $\kow B\subseteq \uqgc\ot B$ relation (\ref{lcoprod}) implies
    \begin{align*}
      \sum_{i,j} \sigma(l^-(c^i_m))J^j_i l^+(c^n_j) - J_m^n\in B^+,\qquad 
      \mbox{for all $m,n=1,\dots,N$},
    \end{align*}
    where $B^+=B\cap \ker(\vep)$.
    By the $\adr(B)$-invariance of $C$ this implies
    \begin{align*}
      \adr \left( \sum_i J^n_i \sigma(l^-(c^i_m))- \sigma^{-1}(l^+(c^n_i))J^i_m 
           \right) (C)=0.
    \end{align*}
    Using relation (\ref{adrcompatible}) one obtains
    \begin{align*}
      \sum_{a,b,i,j,k,l}\left[J^n_i\rf(c^i_l,\sigma(c^b_j))J^j_k
       \rf(c^l_m,c^k_a)-
       \rf(c^b_j,c^n_l)J^j_k\rf(c^k_a,\sigma(c^l_i))J^i_m \right]l(c^a_b)=0
    \end{align*}
    for all $m,n=1,\dots,N$. As $l$ is injective and the $c^a_b$ are 
    linearly independent this yields the desired formula.
    $\blacksquare$
    
    \hspace{.5cm}
    
    \noindent {\bf Remark:}
    Formula (\ref{refl}) coincides with the reflection equation 
    \cite[Equation (5.5)]{a-Dijk96}
    up to replacing $\rtil$ by $R^{-1}$. In Section \ref{leftright}
    we will see how to obtain \cite[Equation (5.5)]{a-Dijk96} by a suitable translation
    from left to right coideal subalgebras.
    
  \subsection{Quantum symmetric pairs}\label{qsp}
    We now wish to apply the above observation to the family of quantum symmetric pairs 
    constructed and investigated by Gail Letzter in a series of papers 
    \cite{a-Letzter99a}, \cite{MSRI-Letzter}, \cite{a-Letzter03}. Hence we first recall the 
    construction and some properties of quantum symmetric pairs.
    
    Let $\theta:\gfrak\rightarrow \gfrak$ be an involutive Lie algebra automorphism,
    i.e.~$\theta^2=\id$, which is maximally split with respect to $\hfrak$ and the 
    chosen triangular decomposition of $\gfrak$ in the sense of 
    \cite[Section 7]{MSRI-Letzter}. 
    Then $\theta$ induces an involution $\Theta$ of the root 
    system $\Delta$ of $\gfrak$ and thus an automorphism of $\hfrak^\ast$. 
    Define $\pt=\{\alpha_i\in \wurz\,|\,\Theta(\alpha_i)=\alpha_i\}$.
    Let $p$ be the permutation of $\{1,\dots,n\}$ such that
       \begin{align}\label{pdefn}
         \Theta(\alpha_i)\in -\alpha_{p(i)}+\Z\pi_{\Theta}\ \mbox{for all}\ \alpha_i\notin \pt
       \end{align}  
    and $p(i)=i$ if $\alpha_i\in \pt$.
    
    Let $\cM\subseteq \uqgc$ denote the Hopf subalgebra generated by 
    $x_i$, $y_i$, $t_i^{\pm 1}$ for $\alpha_i\in \pt$. Define a multiplicative subgroup 
    $\Ttcp\subseteq \uqgc$ by
    \begin{align*}
      \Ttcp=\{\tau(\lambda)\,|\,\lambda \in P(\pi)\,\mbox{and}\,           
                         \Theta(\lambda)=\lambda \}.
    \end{align*}
    By definition the subalgebra $\Bc\subseteq \uqgc$ is generated by $\cM$, $\Ttcp$, and 
    elements $B_i$ for $\alpha_i\in \pi\setminus\pt$ defined by 
    \begin{align}\label{Bidef}
      B_i=y_it_i +d_i\, \thetatil(y_i)t_i +s_i\,t_i
    \end{align}
    for suitable $d_i,s_i\in \cC$, $d_i\neq 0$, and where
    \begin{align*}
      \thetatil(y_i)=(\adr x_{i_m})\dots(\adr x_{i_1})(t_{p(i)}^{-1}x_{p(i)})
    \end{align*}
    for suitable $\alpha_{i_1},\dots,\alpha_{i_m}\in \pt$. 
    For more details, in particular for
    the choice of $d_i$, $s_i$, and $\alpha_{i_1},\dots,\alpha_{i_m}$ consult 
    \cite[Section 7]{MSRI-Letzter}. For the purpose of the present paper it suffices to know
    that $\thetatil(y_i)$ is $(\adr \cM^+)$-invariant. 
    
    Set $B_i=y_it_i$ for $\alpha_i\in \pt$ and set $\cM^+=\cM\cap \uqnp$. Recall the 
    definition of  $\cI$ given in Section \ref{uqgrqg}. Given a multiindex 
    $I=(i_1,\dots, i_m)$ in $\cI$, set $B_I=B_{i_1}\cdots B_{i_m}$.
    By \cite[Section 7]{MSRI-Letzter}, we have a direct sum decomposition	    
    \begin{align}\label{Bdirectsum}
      \Bc=\oplus_{I\in \cI}B_I\cM^+\Ttcp.
    \end{align}

    It is sometimes convenient to replace the generators 
    $\{B_i\,|\,\alpha_i\in \pi\setminus \pt\}$ by a different set of generators
    $\{C_i\,|\,\alpha_i\in \pi\setminus \pt\}$ obtained by interchanging $y_it_i$ and $x_i$.
    More precisely, the new set of generators is obtained via the following map.
    By \cite[Theorem 3.1]{a-Letzter03}, \cite[(1.22)]{a-Letzter-memoirs} there exists
    an $\Cc$-linear algebra antiautomorphism $\kb$ of $\uqgc$ defined by 
    \begin{align*}
      \kb(x_i)=c_{B_i} y_it_i,\quad \kb(y_i)=c_{B_i}^{-1} t_i^{-1}x_i,\quad 
      \kb(\tau(\lambda))=\tau(\lambda)
    \end{align*}
    for suitable nonzero $c_{B_i}\in \Cc$, such that
    \begin{align*}
      \kb(\Bc)=\Bc.
    \end{align*}
    Note that in the above references $\kb$ is defined to be conjugate linear in order 
    to obtain a $\ast$-structure on $\Bc$. However, linearity suffices for the present 
    purposes. Note moreover, that by definition the antipode and the coproduct of $\uqgc$
    satisfy
    \begin{align}\label{kb-sigma-kow}
      \kb(\sigma(u))=\sigma^{-1}(\kb(u)),\quad \kow(\kb(u))=(\kb\ot\kb)\kow(u)
    \end{align}
    for all $u\in \uqgc$. In particular, one gets
    \begin{align}\label{adU-kb}
      \kb((\adr u)C)=(\adr(\sigma^{-1}\circ\kb(u)))\kb(C)
    \end{align}
    for all $u,C\in  \uqgc$, and hence
    \begin{align}\label{adUtau-kb-inv}
      \kb\left((\adr \uqgc)\tau(2\mu)\right)=(\adr \uqgc)\tau(2\mu)
    \end{align}
    for all $\mu\in \Ppiplus$. 
    
    The subalgebra generated by $\cM$ and $\Ttcp$ of $\Bc$ is invariant under $\kb$. Hence
    $\Bc$ is also generated by $\cM$, $\Ttcp$, and elements $C_i=c_{B_i}\kb(B_i)$ for 
    $\alpha_i\in \pi\setminus \pt$. Note that (\ref{adU-kb}) implies
    \begin{align*}
      C_i&=x_i+c_{B_i}d_it_i\kb(\thetatil(y_i))+c_{B_i}s_i t_i\\
         &=x_i+(-1)^m d_i'(\adr y_{i_1})\cdots (\adr y_{i_m})y_{p(i)}+s_i't_i
    \end{align*}
    for suitable $d_i', c_i'\in \Cc$, $d_i'\neq 0$. Formula (\ref{adU-kb}) also implies
    that $\kb(\thetatil(y_i))$ is $(\adr \cM^-)$-invariant because $\thetatil(y_i)$ is 
    $(\adr \cM^+)$-invariant.
    
    %\cite[(3.1), (3.3)]{a-Letzter03} one sees that $B$ is also generated by $\cM$, $\Ttcp$ and
    %elements $C_i$ for $\alpha_i\in \pi\setminus\pt$ defined by 
    %\begin{align}\label{Cidef}
    %      C_i=x_i +q_i^2 d_{p(i)}^{-1}\,\kappa(\thetatil(y_i)) +q_i^2 s_i\,t_i
    %\end{align}
    %where
    %\begin{align*}
    %      \kappa(\thetatil(y_i))=(-1)^m(\adr y_{i_1})\dots(\adr y_{i_m})(y_{p(i)})
    %\end{align*}
    %is $(\adr \cM^-)$-invariant. 
     
    \subsection{The center of quantum symmetric pair coideal subalgebras}\label{cqsp}
    In order to apply Proposition \ref{reflprop} to quantum symmetric pairs we recall some
    results of \cite{a-KLp} on the center $Z(\Bc)$ of the coideal subalgebras 
    $\Bc\subseteq\uqgc$.
    
    Let $W(\pt)\subseteq W$ denote the subgroup generated by the simple reflections
    corresponding to the simple roots in $\pt$. Note that $W(\pt)$ is the Weyl 
    group corresponding to the root system generated by $\pt$. 
    Let $w_0'\in W(\pt)$ denote the longest element in $W(\pt)$.
    As shown in \cite[Theorem 8.5]{a-KLp} the center $Z(\Bc)$ of $\Bc$ has a basis 
    $\{d_\mu\,|\,\mu\in \Ptheta\}$ indexed by the set
    \begin{align*}
       \Ptheta=\{\mu\in \Ppiplus\,|\,\Theta(\mu)=\mu+w_0\mu-w_0'\mu\}.
    \end{align*}
    The set $\Ptheta$ is explicitly determined in \cite[Proposition 9.1]{a-KLp}.
    Note that $\Ptheta=-w_0\Ptheta$.
    By \cite[Theorem 8.3]{a-KLp} the central elements $d_\mu$ for $\mu\in\Ptheta$ satisfy
    \begin{align}\label{dmu}
      d_\mu\in Y^\mu + \sum_{\gamma <2\mutil}\uqgc_{-\gamma}
    \end{align}
    where $Y^\mu\in(\adr \uqgc)\tau(2\mu)$ is a highest weight vector with respect to the
    right adjoint action of spherical weight $2\mutil=\mu-\Theta(\mu)$.
    
    The center $Z(\Bc)$ is invariant under the algebra antiautomorphism $\kb$.
    Recall also from \cite[Theorem 8.3]{a-KLp} that
    \begin{align}\label{dmu=vmu+}
      d_\mu\in v_\mu+\sum_{\nu <\mu} (\adr \uqgc)\tau(2\nu)
    \end{align}
    for some nonzero $v_\mu\in (\adr \uqgc)\tau(2\mu)$ which is uniquely determined 
    up to scalar multiplication. From (\ref{adUtau-kb-inv}) and the fact that
    $\{d_\mu\,|\,\mu\in \Ptheta\}$ is a basis of $Z(\Bc)$ one obtains
    \begin{align}\label{KL8.3}
      \kb(d_\mu)\in a_\mu v_\mu + \sum_{\nu <\mu} (\adr \uqgc)\tau(2\nu)
    \end{align}
    for some $a_\mu\in \Cc\setminus \{0\}$. As $\kb^2=\id$ one even gets $a_\mu\in \{-1,1\}$.
    Hence, replacing $d_\mu$ by $(d_\mu+a_\mu\kb(d_\mu))/2$ we may assume that
    \begin{align}\label{dmu-kb-inv}
      \kb(d_\mu)=a_\mu d_\mu.
    \end{align}
    Define $X^\mu=a_\mu\kb(Y^\mu)$. Then by (\ref{adU-kb}) the element 
    $X^\mu\in (\adr \uqgc)\tau(2\mu)$ is a lowest weight vector of weight $-2\mutil$ 
    with respect to the right adjoint action. 
    Hence (\ref{dmu}) and (\ref{dmu-kb-inv}) imply 
    \begin{align}\label{dmuX}
          d_\mu\in X^\mu + 
          \sum_{\gamma <2\mutil}\uqgc_{\gamma}.
    \end{align}

    \subsection{Minimal weights in $\Ptheta$}\label{minweights}
     By Proposition \ref{reflprop} one can find solutions of the reflection equation via
     elements in $Z(\Bc)\cap(\adr \uqgc)\tau(-2w_0\mu)$. The following notion of minimality
     will lead to the desired central elements. Note that this definition is slightly more
     general than the definition of minimality in say \cite[Exercise 13.13]{b-Humphreys}.
     \begin{definition}
        Let $\gfrak$ be simple. A nonzero weight $\mu\in \Ppiplus$ is called minimal if it 
        satisfies the following property:
        If $\nu\in\Ppiplus$ and $\nu<\mu$ then $\nu=0$.
     \end{definition}
     Note that $\mu\in \Ppiplus$ is minimal if and only if $-w_0\mu$ is minimal. 
     As shown in \cite[Corollary 8.4]{a-KLp} it follows immediately from (\ref{dmu=vmu+}) 
     that for minimal $\mu \in \Ptheta$ we may assume
     \begin{align}\label{dmuin}
        d_\mu\in (\adr \uqgc)\tau(2\mu).
     \end{align}
     In Tables \ref{min} and \ref{minPtheta} we have listed all minimal weights and all 
     minimal weights in $\Ptheta$ using the labeling of simple roots of 
     \cite[11.4]{b-Humphreys}. Symmetric pairs are labeled as in \cite[5.11]{a-Araki62}.
     Following \cite[Section 7]{a-Letzter03}, however, the
     parameter $p$ occurring in Araki's list will be denoted by $r$, 
     and as before $n=\rk(\gfrak)$. 
     Table \ref{min} is obtained by inspection from the table in 
     \cite[p.~69]{b-Humphreys} and Table \ref{minPtheta}
     then follows from \cite[Proposition~9.1]{a-KLp}. 
     \begin{table}[!htb]
       \begin{center}
         \begin{tabular}{c|c||c|c}
          $A_n$ & $\omega_1,\dots,\omega_n, \omega_1{+}\omega_n$
                & $E_7$ & $\omega_1,\,\omega_2$ \\
          $B_n$ & $\omega_1,\omega_n$
                & $E_8$ & $\omega_8$\\
          $C_n$ & $\omega_1,\,\omega_2$
                & $F_4$ & $\omega_4$\\
          $D_n$ & $\omega_1,\,\omega_2,\,\omega_{n-1},\,\omega_n$
                & $G_2$ & $\omega_1$\\ 
	  $E_6$ & $\omega_1,\,\omega_2,\,\omega_6$&&\\
         \end{tabular}
       \end{center}  
      \caption{Minimal weights in $\Ppiplus$}\label{min}
      \end{table}
    
    \begin{table}[!htb]
         \begin{center}
           \begin{tabular}{c|c}
            $AI$, $AII$   &  $\omega_1+\omega_n$ \\
            $AIII$, $AIV$ &  $\omega_1,\dots,\omega_n, \omega_1{+}\omega_n$\\
            $B_n$ & $\omega_1,\omega_n$\\
            $C_n$ & $\omega_1,\,\omega_2$\\
            ($DI$, case 1), $r$ even & $\omega_1,\,\omega_2,\,\omega_{n-1},\,\omega_n$\\
            ($DI$, case 1), $DII$, $r$ odd & $\omega_1,\,\omega_2$\\
            ($DI$, case 2,3), ($DIII$, case 1), $n$ even & $\omega_1,\,\omega_2,\,\omega_{n-1},\,\omega_n$\\
            ($DI$, case 2,3), ($DIII$, case 1), $n$ odd & $\omega_1,\,\omega_2$\\
            ($DIII$, case 2) & $\omega_1,\,\omega_2, \,\omega_{n-1},\,\omega_n$\\
    	      $EI$, $EIV$ & $\omega_2$\\
    	      $EII$, $EIII$ & $\omega_1,\,\omega_2,\,\omega_6$\\
    	      $E_7$ & $\omega_1,\,\omega_2$ \\
    	      $E_8$ & $\omega_8$\\
    	      $F_4$ & $\omega_4$\\
    	      $G_2$ & $\omega_1$
             \end{tabular}
           \end{center}  
          \caption{Minimal weights in $\Ptheta$}\label{minPtheta}
      \end{table}

  \subsection{Invariants} \label{invariants} 
    Assume now that $\mu\in \Ptheta$ and let $L_\mu\subseteq \uqgc$ denote the left 
    coideal generated by the central element $d_{-w_0\mu}$. As
    $\Bc\subseteq \uqgc$ is a left coideal one has $L_\mu\subseteq \Bc$.
    If $\mu\in \Ptheta$ is minimal
    then (\ref{dmuin}) holds also with $\mu$ replaced by $-w_0\mu$. 
    In this case let again $J$ denote
    the $N\times N$-matrix defined by $d_\mu=\sum_{i,j}J^j_il(c^i_j)$.
    By Proposition \ref{reflprop} the matrix $J$ is a solution of the reflection equation
    (\ref{refl}).
    Let $L^+$ and $L^-$ denote the $N\times N$ matrix with entries $l^+(c^i_j)$ and
    $l^-(c^i_j)$, respectively. By construction $L_\mu$ is the linear span of the 
    matrix entries of  $\sigma(L^-)^tJ^t(L^+)^t$. 
    Up to conventional changes explained in Section \ref{leftright} the coideal $L_\mu$
    is of the type considered by M.~Noumi, T.~Sugitani, and M.~Dijkhuizen.   
    
    For any left $\uqgc$-module $M$ and any subset $L\subseteq \uqgc$ define 
    $M^L=\{m\in M\,|\,am=\vep(a)m\,\forall a\in L\}$. We conjecture that
    $L_\mu$ is big enough in the following sense.
    \begin{conjecture}\label{conj1}
      Assume that $\gfrak$ is simple and  that $\mu\in \Ptheta$ is minimal. 
      Then $V(\nu)^\Bc=V(\nu)^{L_\mu}$ for all $\nu\in\Ppiplus$.
    \end{conjecture}
    For the classical symmetric pairs and a suitable choice of a minimal $\mu\in \Ppiplus$
    the coideals $L_\mu$ have been explicitly constructed. In this case the claim of the
    conjecture holds by \cite[Theorem 1]{a-NS95} and \cite[Theorem 7.7]{MSRI-Letzter}.
    In Section \ref{towards} we will obtain general results aiming at a proof of 
    Conjecture \ref{conj1}. Although the general case remains open these results are strong
    enough to verify Conjecture \ref{conj1} for the symmetric pair $FII$.
    
    To obtain solutions of the reflection equation via Proposition \ref{reflprop} it is
    necessary that $d_\mu$ satisfies (\ref{dmuin}) and hence it is natural to assume
    that $\mu$ is minimal. However, one could ask if Conjecture \ref{conj1} also holds
    for not necessarily minimal $\mu\in \Ptheta$. This also allows us to consider
    irreducible symmetric pairs $(\gfrak,\gfrak^\theta)$ with not necessarily simple 
    $\gfrak$.
    \begin{conjecture}\label{conj2}
      Assume that $(\gfrak, \gfrak^\theta)$ is irreducible and that $\mu\in \Ptheta$. 
      Then $V(\nu)^\Bc=V(\nu)^{L_\mu}$ for all $\nu\in\Ppiplus$.
    \end{conjecture}
    In Section \ref{pt-empty} the claim of Conjecture \ref{conj2} is verified for all quantum 
    symmetric pairs with $\pt=\emptyset$. In particular, Conjecture \ref{conj2} is seen 
    to hold in the case where $\gfrak=\gfrak'\oplus\gfrak'$ for some simple $\gfrak'$
    and $\theta$ interchanges the two isomorphic components.
    
    \medskip
    
    \noindent {\bf Remark:} 
    M.~Noumi constructed quantum symmetric pairs of type AI and AII using the $R$-matrix
    in the vector representation $V(\omega_1)$ and a reflection equation
    \cite[(2.7)]{a-Noumi96} different from (\ref{refl}) and (\ref{Dijkrefl}). 
    However, $\omega_1\notin P_{Z(\Bc)}$ for AI and AII. The construction 
    in \ref{Coideal-Refl} only yields a solution of the reflection equation for the 
    $R$-matrix of $V(\omega_1+\omega_n)$. Conjectures \ref{conj1} and \ref{conj2} claim
    that with respect to its invariants the coideal constructed here is equivalent to
    the coideal constructed by Noumi. Noumi's reflection equation
    occurs for all classical symmetric pairs apart from AIII \cite{a-NS95}.
    However, for $\gfrak$ of type $B$, $C$, or $D$
    the vector representation $V(\omega_1)$ is self dual and hence Noumi's 
    reflection equation can be translated into (\ref{refl}) in these cases.
    
    In view of the above remarks it is one of the conceptual insights of the present paper
    that Noumi-Sugitani-Dijkhuizen type constructions can be performed with one and
    the same reflection equation uniformly for all symmetric pairs.

    \subsection{First properties of $L_\mu$}
    The following two observations will be frequently used in Sections 
    \ref{towards} and \ref{pt-empty}.
    \begin{lemma}\label{adlM-inv}
      For any $\mu\in \Ptheta$ the left coideal $L_\mu$ is $(\adl \cM)$-invariant.
    \end{lemma}
    {\bf Proof:} For $b\in \cM\cap \ker(\vep)$ and $C\in Z(\Bc)$ one has analogously to
    (\ref{kowadr})
       \begin{align*}
         0=\kow((\adl b)C)=b_{(1)}C_{(1)}\sigma(b_{(3)})\ot (\adl b_{(2)}) C_{(2)}.
       \end{align*}
    In particular for $b=x_i$, $\alpha_i\in \pt$, and $C=d_{-w_0\mu}$ one obtains
      \begin{align*}
         C_{(1)}\ot (\adl x_i)C_{(2)}=t_i^{-1}(x_iC_{(1)}-C_{(1)}x_i)\ot C_{(2)}\in 
         \uqgc\ot L_\mu.
      \end{align*}
    Hence $L_\mu$ is $(\adl x_i)$-invariant. Similarly one obtains that $L_\mu$ is 
    $(\adl y_i)$-invariant and $(\adl t_i)$-invariant if $\alpha_i\in \pt$.\, $\blacksquare$
    
    \begin{lemma}\label{kb-inv}
       For any $\mu\in \Ptheta$ the left coideal $L_\mu$ is invariant under the 
       algebra antiautomorphism $\kb$.
    \end{lemma}
    {\bf Proof:} This follows from (\ref{dmu-kb-inv}) and the second formula of
                 (\ref{kb-sigma-kow}). $\blacksquare$

  \subsection{Left versus right}\label{leftright}
    In this section the result from \ref{Coideal-Refl} is translated into the 
    setting of right coideal subalgebras. This reproduces the reflection equation
    \cite[(5.5)]{a-Dijk96}. The results of the present section are not relevant
    for the rest of the paper.
  
    Using the antipode of $\uqgc$ there are two ways to turn a left
    $\uqgc$-module $V$ into a right $\uqgc$-module. One can define
    a right action on $V$ by
    \begin{align}\label{v+action}
      v\ract_+u:=\sigma(u)v\qquad \forall\, u\in\uqgc,\, v\in V
    \end{align}
    or by
    \begin{align}\label{v-action}
          v\ract_-u:=\sigma^{-1}(u)v\qquad \forall\, u\in\uqgc,\, v\in V.
    \end{align}
    We will denote the vectorspace $V$ endowed with the right module 
    structures (\ref{v+action}) and (\ref{v-action}) by $V_r^+$ 
    and $V_r^-$, respectively.
    For finite dimensional $V$ the right $\uqgc$-modules $V^+_r$ and 
    $V^-_r$ are isomorphic, however, the isomorphism is not given 
    by the identity map of the underlying vectorspace $V$.
    
    This little observation explains why the reflection equation
    (\ref{refl}) slightly differs from the reflection equation 
    \cite[(5.5)]{a-Dijk96}. To be more precise, note that 
    by (\ref{l-compatible}) the map
    \begin{align}\label{-+iso}
      V(\lambda)^\ast\ot V(\lambda)^+_r &\rightarrow 
         (\adr \uqgc)(\tau(2\lambda)),
         \qquad v^\ast_i\ot v_j\mapsto l(c^i_j)
    \end{align}
    is an isomorphism of right $\uqgc$-modules.
    For any $N\times N$-matrix $J$ let $B_J\subseteq \uqgc$ denote the 
    left coideal subalgebra generated by the matrix entries of 
    \begin{align*}
      \sigma(L^-)^tJ^t(L^+)^t.
    \end{align*}
    With this notation the proof of Proposition \ref{reflprop} also 
    proves the following generalization of a right version of 
    \cite[Prop.~2.3]{a-Noumi96}, \cite[Prop.~6.5]{a-DS99}.
    \begin{proposition}\label{invreflprop}
      For any $N\times N$-matrix $J$ the element 
      \begin{align*}
        \sum_{i,j}J^j_i\, v_i^\ast\ot v_j\in 
        V(\lambda)^\ast\ot V(\lambda)^+_r
      \end{align*}
      is $B_J$-invariant if and only if $J$ satisfies the 
      reflection equation (\ref{refl}).
    \end{proposition}
    {\bf Proof:} In view of the isomorphism (\ref{-+iso}) all
    that needs to be checked is that $\sum_{i,j}J^j_i l(c^i_j)$
    is $(\adr B_J)$-invariant if and only if $J$ satisfies 
    (\ref{refl}). This has indeed been verified in the proof of 
    Proposition \ref{reflprop}.
    $\blacksquare$
    
    \vspace{.5cm}
    
    To obtain the reflection equation from \cite[(5.5)]{a-Dijk96},
    \cite[(6.4)]{a-DS99} one has to write the invariant element as an 
    element in $V(\lambda)^\ast\ot V(\lambda)^-_r$.
    Note that there is an isomorphism of right $\uqgc$-modules given by
    \begin{align}
      \eta_{-+}: V(\lambda)_r^-\rightarrow V(\lambda)^+_r,\qquad
      v_\nu \mapsto q^{(\lambda-\nu,2\rho)} v_\nu
    \end{align}
    where $v_\nu\in V(\lambda)_\nu$ denotes a weight vector of weight
    $\nu$ in the left $\uqgc$-module $V(\lambda)$ and 
    $\rho=\frac{1}{2}\sum_{\alpha\in \Delta^+}\alpha$.
    Indeed, by (\ref{squareofAnti}) for any homogeneous element
    $a\in \uqgc_\mu$ on has
    \begin{align*}
      v_\lambda\ract_- a=\sigma^{-1}(a)v_\lambda=
      q^{(\mu,2\rho)}\sigma(a)v_\lambda
      =q^{(\mu,2\rho)}v_\lambda\ract_+ a
    \end{align*}
    and therefore
    \begin{align*}
      \eta_{-+}(v_\lambda \ract_- a)
       =q^{(-\mu,2\rho)}\sigma^{-1}(a)v_\lambda
       =v_\lambda\ract_+ a.
    \end{align*}
    For any $N\times N$-matrix $\Jtil$ define $\Btil_\Jtil:=B_J$ where
    $J^j_i:=q^{(\lambda-\wght(v_j),2\rho)}\Jtil^j_i$.
    The following proposition states that exchange of left and right
    by means of the antipode $\sigma$ transforms the situation
    considered here precisely into the situation considered in  
    \cite{a-DS99} for the special case of Grassmannians. 
    \begin{proposition}
      For any $N\times N$-matrix $\Jtil$ the element
      \begin{align*}
              C:=\sum_{i,j}\Jtil^j_i\, v_i^\ast\ot v_j\in 
              V(\lambda)^\ast\ot V(\lambda)^-_r
      \end{align*}
      is $\Btil_\Jtil$-invariant if and only if $\Jtil$ satisfies the
      reflection equation
       \begin{align}
	      \Jtil_1 R^{-1} \Jtil_2 R = 
	       R_{21} \Jtil_2 R^{-1}_{21} \Jtil_1. \label{Dijkrefl}
       \end{align} 
       Moreover, let $B_\Jtil^D\subseteq\uqgc$ denote the right
       coideal subalgebra generated by the coefficients of the 
       $N\times N$-matrix $\sigma(L^+)\Jtil L^-$. Then
       $\Btil_\Jtil=\sigma(B^D_\Jtil)$.
    \end{proposition}
    {\bf Proof:}
    Note that
    \begin{align*}
      (\id\ot\eta_{-+})\left(\sum_{i,j}\Jtil^j_i\, v^\ast_i\ot v_j\right)=
      \sum_{i,j}J^j_i\, v^\ast_i\ot v_j
      \in V(\lambda)^\ast\ot V(\lambda)^+_r
    \end{align*}
    where $J^j_i=q^{(\lambda-\wght(v_j),2\rho)}\Jtil^j_i$. Thus 
    by Proposition \ref{invreflprop} the element $C$ is 
    $\Btil_\Jtil$-invariant if and only if $J$ satisfies the reflection equation
    (\ref{refl}). Thus $C$ is $\Btil_\Jtil$-invariant if and only if 
    $\Jtil$ satisfies
    \begin{align}
         \sum_{i,j,k,l}&q^{(\lambda-\wght v_n,2\rho)}\Jtil^n_i
          \rf(c^i_l,\sigma(c^b_j))q^{(\lambda-\wght v_j,2\rho)}
           \Jtil^j_k \rf(c^l_m,c^k_a)\nonumber\\
           &-\rf(c^b_j,c^n_l)q^{(\lambda-\wght v_j,2\rho)}\Jtil^j_k
           \rf(c^k_a,\sigma(c^l_i))q^{(\lambda-\wght v_i,2\rho)}
           \Jtil^i_m =0 \label{refl1}
    \end{align}
    for all $a,b,m,n$.
    Note that (\ref{squareofAnti}) implies 
    \begin{align}\label{sigmasquare}
      q^{(\lambda-\wght v_j,2\rho)}c^i_j
        =q^{(\lambda-\wght v_i,2\rho)}\sigma^{-2}(c^i_j).
    \end{align}  
    Using (\ref{sigmasquare}) relation (\ref{refl1}) is seen to be
    equivalent to
    \begin{align*}
              \sum_{i,j,k,l}&q^{(\lambda-\wght v_n,2\rho)}
              q^{(\lambda-\wght v_b,2\rho)}\Jtil^n_i
              \rf(c^i_l,\sigma^{-1}(c^b_j))
               \Jtil^j_k \rf(c^l_m,c^k_a)\\
               &-q^{(\lambda-\wght v_n,2\rho)}
              q^{(\lambda-\wght v_b,2\rho)}
              \rf(\sigma^{-2}(c^b_j),\sigma^{-2}(c^n_l))\Jtil^j_k
               \rf(c^k_a,\sigma^{-1}(c^l_i))\Jtil^i_m =0
    \end{align*}
    for all $a,b,m,n$. Using Proposition \ref{r-trig} iv)
    one verifies that the above equation is equivalent to
    \begin{align*}
                  \sum_{i,j,k,l}\Jtil^n_i
                  \rf(\sigma(c^i_l),(c^b_j))
                   \Jtil^j_k \rf(c^l_m,c^k_a)
                   -\rf(c^b_j,c^n_l)\Jtil^j_k
                   \rf(\sigma(c^k_a),c^l_i)\Jtil^i_m =0
    \end{align*}
    for all $a,b,m,n$. This is just the desired equation
    (\ref{Dijkrefl}) written explicitly.
    
    To verify the second statement, note that (\ref{sigmasquare})
    implies
    \begin{align*}
      \sum_{i,j}\sigma(l^-(c^i_m)) J^j_i\, l^+(c^n_j)&=
      \sum_{i,j}\sigma(l^-(c^i_m))\Jtil^j_i\, \sigma^2(l^+(c^n_j))
      q^{(\lambda-\wght v_n,2\rho)}\\
      &=\sigma\left(\sum_{i,j}\sigma(l^+(c^n_j))\Jtil^j_i\, l^-(c^i_m)
      \right)q^{(\lambda-\wght v_n,2\rho)}.\quad\blacksquare
    \end{align*}
  
 \section{Towards a proof of Conjecture \ref{conj1}}\label{towards}   
   Throughout this section we assume the $\gfrak$ is simple. The proof of the stronger
   Conjecture \ref{conj2} in the case $\gfrak=\gfrak'\oplus\gfrak'$  is contained in 
   Section \ref{pt-empty}.
   
   \subsection{Coproduct of central elements}
     Recall the definition of the projection maps $\pi_{\alpha,\beta}$ and $\pi_\alpha$ for 
     $\alpha,\beta\in Q(\pi)$ from Section \ref{uqgrqg}. By the following Lemma the coproduct of
     the central element $d_\mu$ for $\mu\in \Ptheta$ can be in part read off the element
     $Y^\mu$ defined by (\ref{dmu}).

    \begin{lemma}\label{dYcoprod}
       Let $\mu\in \Ptheta$ and $\alpha,\beta\in Q(\pi)$. Then
       \begin{align*}
         (\pi_{-(\mu-w_0\mu-\alpha),\beta}\ot \mbox{\upshape \id})\kow (d_\mu)\in &
         (\pi_{-(\mu-w_0\mu-\alpha),\beta}\ot \pi_{\mu-w_0'\mu-\alpha-\beta})\kow (Y^\mu)\\
         &+ \sum_{\eta<w_0'\mu-\mu+\alpha+\beta}\uqnm_{-(\mu-w_0\mu-\alpha)}\Uc^0 G^+_\beta \ot
         \uqgc_{-\eta}.
       \end{align*}
   \end{lemma}
   \noindent {\bf Proof:} Recall from (\ref{dmu}) that 
   $d_\mu\in Y^\mu+\sum_{\gamma<2\mutil}\uqgc_{-\gamma}$ 
   and $Y^\mu\in \uqgc_{-2\mutil}$. Note moreover that 
   $\kow \uqgc_\nu\subseteq \sum_{\zeta+\xi=\nu} \uqgc_\zeta\ot \uqgc_\xi$. Hence using
   $2\mutil=w_0'\mu-w_0\mu$ one obtains
   \begin{align*}
     (\pi_{-(\mu-w_0\mu-\alpha),\beta}\ot \mbox{\upshape \id})\kow (Y^\mu) &=
     (\pi_{-(\mu-w_0\mu-\alpha),\beta}\ot \pi_{\mu-w_0'\mu-\alpha-\beta})\kow (Y^\mu)\\
     (\pi_{-(\mu-w_0\mu-\alpha),\beta}\ot \mbox{\upshape \id}) \kow (\uqgc_{-\gamma}) &\subseteq
     \sum_{\eta<w_0'\mu-\mu+\alpha+\beta}\uqnm_{-(\mu-w_0\mu-\alpha)}\Uc^0 G^+_\beta \ot
     \uqgc_{-\eta}
   \end{align*}
   for any $\gamma<2\mutil$.
   $\blacksquare$

   \subsection{Finding one generator $B_i$ in $L_{\mu}$}
   Recall that for any $\mu\in \Ptheta$ the symbol $L_\mu$ denotes the left coideal generated 
   by the central element $d_{-w_0\mu}\in Z(\Bc)$.
   
   \begin{proposition}\label{BiinLmu}
     Let $\mu\in \Ptheta$ such that $(\mu,\alpha_j)=0$ for all $\alpha_j\in \pt$ and
     $2(\mu,\alpha_i)/(\alpha_i,\alpha_i)=1$ for some $\alpha_i\notin \pt$. 
     Then $B_i\tau({-}\Theta(\mu){-}\mu)\in L_\mu$ and 
     $C_i\tau({-}\Theta(\mu){-}\mu)\in L_\mu$.
   \end{proposition}
   
   \noindent {\bf Proof:} 
   Consider the central element $d_{-w_0\mu}$ which lies in $(\adr \uqgc)\tau(-2w_0\mu)$
   up to terms of lower filter degree, and its highest weight component $Y^{-w_0\mu}$.
   %Note that 
   %\begin{align}\label{considerY}
   %  (\pi_{\mu-w_0\mu-\alpha_i,0}\ot \pi_{-\alpha_i})\kow (d_{-w_0\mu}) 
   %  =(\pi_{\mu-w_0\mu-\alpha_i,0}\ot \pi_{-\alpha_i})\kow (Y^{-w_0\mu}). 
   %\end{align}  
   Note that $(\mu,\alpha_j)=0$ for all $\alpha_j\in \pt$ implies $w_0'\mu=\mu$. 
   Hence, using $\mu\in \Ptheta$, one obtains $\Theta(w_0\mu)=\mu$ and thus
   \begin{align}\label{tildew0mu} 
     -2\widetilde{w_0\mu}=-w_0\mu-\Theta(-w_0\mu)=\mu-w_0\mu.
   \end{align}  
   By \cite[4.8]{a-JoLet2} there is an isomorphism of right $\uqnm$-modules
   \begin{align}\label{U-iso}
     (\adr \uqnm)\tau(-2w_0\mu)\cong V(\mu)^\ast, \quad (\adr u)\tau(-2w_0\mu)
                         \mapsto v_{w_0\mu}^\ast u
   \end{align}
   where $u\in \uqnm$ and $v_{w_0\mu}^\ast\in V(\mu)^\ast$ is nonzero on the lowest 
   weight vector $v_{w_0\mu}$ and vanishes on all other weight vectors of $V(\mu)$.
   The isomorphism (\ref{U-iso}) implies that
   \begin{align*}
     Y^{-w_0\mu}=(\adr y_i)(\adr \Ytil)\tau(-2w_0\mu)
   \end{align*}
   for some $\Ytil\in \uqnm_{-(\mu-w_0\mu-\alpha_i)}$. Moreover, using
   $2(\mu,\alpha_i)/(\alpha_i,\alpha_i)=1$
   and the isomorphism (\ref{U-iso}) one sees that the weight space with respect to the left 
   adjoint action $((\adr \uqnm)\tau(-2w_0\mu))_{-(\mu-w_0\mu-2\alpha_i)}$ is empty.
   Hence, by (\ref{adryC}) one obtains
   \begin{align*}
     &(\pi_{-(\mu-w_0\mu-\alpha_i),0}\ot \pi_{-\alpha_i})\kow (Y^{-w_0\mu})\\
      =&(\adr \Ytil)\tau(-2w_0\mu)\ot [\tau({-}w_0\mu{-}\mu{+}\alpha_i)y_i 
         -q^{-(\alpha_i,\mu{-}w_0\mu{-}\alpha_i)} y_i\tau({-}w_0\mu{-}\mu{+}\alpha_i)]\\
      =&(1-q^{(\alpha_i,\alpha_i)})(\adr \Ytil)\tau(-2w_0\mu)\ot 
        \tau({-}w_0\mu{-}\mu{+}\alpha_i)y_i\\
      =& (q^{-(\alpha_i,\alpha_i)}-1) (\adr \Ytil)\tau(-2w_0\mu)\ot 
        \tau({-}\mu{-}\Theta(\mu))y_it_i
   \end{align*}
   Lemma \ref{dYcoprod} for $\alpha=\alpha_i$, $\beta=0$ and the coideal property for 
   $L_\mu$ now imply that
   \begin{align*}
      L_\mu \cap \left(  y_i t_i\tau(-\Theta(\mu)-\mu) +\sum_{\eta<\alpha_i}
     \uqgc_{-\eta}\right)\neq 0.
   \end{align*}
   Hence by (\ref{Bdirectsum}) there exists $m\in \cM\Ttcp$ such that 
   $B_i\tau({-}\Theta(\mu){-}\mu)+m\in L_\mu$.
   Using the coideal property again one obtains $B_i\tau({-}\Theta(\mu){-}\mu)\in L_\mu$.
   The relation $C_i\tau({-}\Theta(\mu){-}\mu)\in L_\mu$ now follows from Lemma \ref{kb-inv}.
   $\blacksquare$

   \subsection{Finding generators $x_j,y_j\in L_\mu\Ttcp\cap \cM$}
   \begin{proposition}\label{xjyj}
        Let $\mu\in\Ptheta$ and $\alpha_i\in \pi\setminus\pt$ and $\tau(\lambda)\in \Ttcp$.
        Assume that $B_i\tau(\lambda)\in L_\mu$ and 
        $(\alpha_i,\alpha_j)\neq 0$ for some $\alpha_j\in \pt$. 
        Then $x_j\tau(\lambda+\alpha_i+\Theta(\alpha_i))\in L_\mu$
        and $y_jt_j\tau(\lambda+\alpha_i+\Theta(\alpha_i))\in L_\mu$.
      \end{proposition}
   {\bf Proof:} Consult the list \cite[5.11]{a-Araki62} of irreducible symmetric pairs 
   for simple $\gfrak$ and note that $(\alpha_i,\alpha_j)\neq 0$ implies
   \begin{align}\label{alpha_irel}
     2(\alpha_i,\alpha_j)/(\alpha_j,\alpha_j)=-1
   \end{align}  
   unless the symmetric pair is of one of the following types:
   \begin{enumerate}
   \item ($CII$, case 2) and $i=n$, $j=n-1$.
   \item  $BI/BII$ and $r=i=n-1$, $j=n$.
   \end{enumerate}
   In case 1.~by \cite[Section 7]{a-Letzter03} one has 
   $B_n=y_n t_n +[(\adr x_{n-1}^2) t_n^{-1} x_n]t_n$. Hence, using
   (\ref{kowadr}) and the relation $(\alpha_n,\alpha_{n-1})=-(\alpha_{n-1},\alpha_{n-1})$ 
   one computes
   \begin{align*}
     (\pi_{\alpha_{n-1}+\alpha_n}\ot \id)\kow B_n=(q^{2(\alpha_{n-1},\alpha_{n-1})}-1)
     [(\adr x_{n-1}) t_n^{-1}x_n] t_n\ot t_{n-1}^{-2}x_{n-1}.
   \end{align*}
   The coideal property of $L_\mu$ now implies
   \begin{align*}
     x_{n-1}\tau(\lambda+\alpha_n+\Theta(\alpha_n))= 
     x_{n-1}\tau(\lambda-2 \alpha_{n-1})\in L_\mu.
   \end{align*}
   The relation $ y_{n-1}t_{n-1}\tau(\lambda+\alpha_n+\Theta(\alpha_n))\in L_\mu$ follows
   from Lemma \ref{kb-inv}. Interchanging $n$ and $n-1$ one obtains the claim of the 
   Lemma in the case 2. Hence from now on we may assume (\ref{alpha_irel}) to hold.
      
   Note that $t_{p(i)}^{-1}x_{p(i)}$ generates a finite dimensional $\adr(\cM)$-module
   $V(\cM,i)$ of lowest weight $\Theta(\alpha_i)$ with respect to the right adjoint action. 
   The relation (\ref{alpha_irel}) implies that 
   \begin{align}\label{weightspacezero}
     V(\cM,i)_{-\Theta(\alpha_i)-2\alpha_j}=\{0\}.
   \end{align}
   and that there exists $\Xtil\in \cM^+$ such that 
   \begin{align*}
     B_i=y_it_i+c\left((\adr x_j) (\adr\Xtil) (t_{p(i)}^{-1}x_{p(i)})  \right)t_i +s_it_i.
   \end{align*}
   Using (\ref{adrxC}), (\ref{weightspacezero}), and the abbreviation 
   $t=\tau(\Theta(\alpha_i)+\alpha_j)$ one obtains
   \begin{align*}
      (\pi_{-\Theta(\alpha_i)-\alpha_j} & \ot \id)(B_i\tau(\lambda))\\
      =&[(\adr \Xtil)t_{p(i)}^{-1}x_{p(i)}]t_i\tau(\lambda)\ot 
        [\sigma(x_j) t-q^{(\alpha_j,\Theta(\alpha_i)+\alpha_j)} 
         t\sigma(x_j)]t_i\tau(\lambda)\\
      =& -[(\adr \Xtil)t_{p(i)}^{-1}x_{p(i)}]t_i\tau(\lambda)\ot 
         t_j^{-1}x_jtt_i\tau(\lambda)(1-q^{2(\alpha_j,\Theta(\alpha_i)+\alpha_j)}).
   \end{align*}
   Note that (\ref{alpha_irel}) implies 
   $(\alpha_j,\Theta(\alpha_i)+\alpha_j)=(\alpha_j,\alpha_i+\alpha_j)=
     -(\alpha_j,\alpha_i)\neq 0$. Hence the coideal property of $L_\mu$ implies
   $x_j\tau(\lambda+\alpha_i+\Theta(\alpha_i))\in L_\mu$.
   The relation 
   $y_jt_j\tau(\lambda+\alpha_i+\Theta(\alpha_i))\in L_\mu$ now follows from Lemma
   \ref{kb-inv}. $\blacksquare$

   \subsection{Obtaining $\cM$-invariance}
   By \cite[(7.5)]{MSRI-Letzter} for any $\alpha_i\in \pi\setminus \pt$ one can write
   \begin{align*}
     \Theta(\alpha_i)+\alpha_{p(i)}=\sum_{\alpha_j\in \pt} n_{ij} \alpha_j
   \end{align*}
   for some uniquely determined nonpositive integers $n_{ij}$. Define a subset of $\pt$ by
   \begin{align*}
     \pi_{\Theta,i}=\{\alpha_j\in \pt\,|\, n_{ij}\neq 0\}
   \end{align*}
   and let $\cM_i$ denote the subalgebra of $\cM$ generated by  
   $\{x_j,y_j,t_j,t_j^{-1}\,|\,\alpha_j{\in} \pi_{\Theta,i}\}$.
   
   \begin{proposition}\label{M-inv}
           Let $\mu\in\Ptheta$ and $\alpha_i\in \pi\setminus\pt$ and $\tau(\lambda)\in \Ttcp$.
           Assume that $B_i\tau(\lambda)\in L_\mu$. Then 
           \begin{align*}
             V(\nu)^{L_\mu}\subseteq V(\nu)^{\cM_i}
           \end{align*}
           for any $\nu\in \Ppiplus$.
   \end{proposition}
   {\bf Proof:} Assume $\alpha_m\in \pti$. There exist uniquely determined pairwise distinct 
   simple roots $\alpha_i=\alpha_{i_0},\alpha_{i_1},\dots, \alpha_{i_k}=\alpha_m$ such that 
   $(\alpha_{i_j},\alpha_{i_{j'}})\neq 0$ if and only if $|j-j'|=1$.
   For any $v\in V(\nu)^{L_\mu}$ we prove $x_m v=0=y_m v$ by induction on $k$. Proposition
   \ref{xjyj} implies $x_{i_1}\tau(\lambda+\alpha_i+\Theta(\alpha_i))\in L_\mu$
        and $y_{i_1}t_{i_1}\tau(\lambda+\alpha_i+\Theta(\alpha_i))\in L_\mu$,
   and hence the claim holds for $k=1$. Moreover, by Lemma \ref{adlM-inv} one has
   \begin{align}\label{adlM}
     \adl(x_{i_k}x_{i_{k-1}}\dots x_{i_2})
              x_{i_1}\tau(\lambda+\alpha_i+\Theta(\alpha_i))\in L_\mu.
   \end{align}
   By induction hypothesis $v$ is invariant under $x_{i_1},\dots,x_{i_{k-1}}$. 
   Hence (\ref{adlM}) implies
   \begin{align}\label{x0}
     x_{i_1} x_{i_2} \dots x_{i_{k-1}} x_{i_k} v=0.
   \end{align}
   Moreover, as the $\alpha_{i_j}$ are pairwise distinct and $v$ is invariant under $y_{i_1}$
   one obtains
   \begin{align}\label{y0}
        y_{i_1} x_{i_1}\dots x_{i_{k-1}} x_{i_k} v=0.
   \end{align}
   Relations (\ref{x0}), (\ref{y0}) and the induction hypothesis imply
   \begin{align*}
     x_{i_2} \dots x_{i_{k-1}} x_{i_k} v= t_{i_1}x_{i_2} \dots x_{i_{k-1}} x_{i_k} v
     =q^{(\alpha_{i_1},\alpha_{i_2})}x_{i_2} \dots x_{i_{k-1}} x_{i_k} v
   \end{align*}
   and hence $x_{i_2} \dots x_{i_{k-1}} x_{i_k} v=0$. Repeating the above argument
   one obtains $x_{i_k}v=0$. The relation $y_{i_k}v=0$ is obtained in a similar manner.
   $\blacksquare$

   \subsection{The example $FII$}
   We apply the results of the previous subsections to prove Conjecture~\ref{conj1} for the 
   quantum symmetric pair of type $FII$. Let $\mu=\omega_4$ be the only minimal weight
   in $\Ptheta$, see Table~\ref{minPtheta}. Note that $\Theta(\omega_4)=-\omega_4$ and
   hence Proposition~\ref{BiinLmu} implies $B_4,C_4\in L_\mu$. By Proposition~\ref{M-inv}
   one obtains the claim of Conjecture~\ref{conj1}. 
   
   Note that we have not proved that $x_i, y_i\in L_\mu \Ttcp$ for $i=1,2$. However, 
   Proposition \ref{xjyj} implies that 
   $x_3, y_3t_3\in L_\mu\tau(\alpha_1+2\alpha_2+3\alpha_3)$.

\medskip
   
   \noindent{\bf Remark:} Note that in Propositions \ref{BiinLmu}, \ref{xjyj}, and \ref{M-inv} we have not
   assumed that $\mu\in \Ptheta$ is minimal. However, the assumption 
   $2(\mu,\alpha_i)=(\alpha_i,\alpha_i)$ of Proposition \ref{BiinLmu} doesn't hold for general
   $\mu\in \Ptheta$. In general the assumption $(\mu,\alpha_j)=0$ for all $\alpha_j\in \pt$ 
   doesn't even hold for all minimal $\mu\in \Ptheta$. This is the reason why our arguments 
   are not yet good enough to prove Conjectures \ref{conj1} and \ref{conj2} in general. 

 \section{The proof of Conjecture \ref{conj2} for $\pt=\emptyset$}\label{pt-empty} 
  Let as before $L_{-w_0\mu}$ denote the left coideal 
  generated by the central element $d_\mu$. 
  For any $\alpha_i\in \pi$ we want to show that there exists $t\in T_\Theta'$
  such that $B_i t\in L_{-w_0\mu}$. For any multiindex $I$ recall the notation $x_I$, $y_I$ from Section
  \ref{uqgrqg}.
 
 \subsection{One dimensional weight spaces of $V(\mu)^\ast$}
 Fix $\mu\in \Ppiplus$ and let $\cIm$ denote a set of multiindices, including the empty 
 multiindex $I_0=()$, such that $\{v_\mu^\ast x_I\,|\,I\in \cIm\}$
 is a basis of $V(\mu)^\ast$. Here $v_\mu^\ast\in V(\mu)^\ast$ is defined by 
 $v_\mu^\ast(v_\mu)=1$ and $v_\mu^\ast(v)=0$ for all $v\in V(\mu)$ of weight unequal $\mu$.
 In short, $v_\mu^\ast$ is a highest weight vector of the right $\uqgc$-module $V(\mu)^\ast$.
 Note that $\{v_{-\mu}^\ast y_I\,|\,I\in \cIm\}$ is a basis of $V(-w_0\mu)^\ast$ where
 $v_{-\mu}^\ast$ denotes a lowest weight vector of $V(-w_0\mu)^\ast$.
 In general, given $\mu\in \Ppiplus$, there are many possible choices for $\cIm$. However,
 certain multiindices always have to be contained, as expressed by the following Lemma.
 
 \begin{lemma}\label{IinI}
   Let $\mu\in \Ppiplus$ and let $I=(i_1,\dots,i_k)$ be a multiindex with the following 
   properties:
   \begin{itemize}
     \item[1)] $(\mu,\alpha_{i_r})\neq 0$ if and only if $r=1$.
     \item[2)] $(\alpha_{i_r},\alpha_{i_s})\neq 0$ if and only if $|r{-}s|=1$.
   \end{itemize}
    Then $I\in \cIm$. 
    Furthermore, the multiindices $I_l:=(i_1,\dots,i_l)$ for $1\le l \le k$ satisfy the 
    relations
   \begin{align}\label{dim1ws1}
     \dim(V(\mu)^\ast_{\mu-\wght(I_l)})=1=\dim(V(-w_0\mu)^\ast_{-\mu+\wght(I_l)}).
   \end{align}  
   If moreover $\mu\in \Ptheta$ then
   \begin{align}\label{dim1ws2}
        \dim(V(\mu)^\ast_{\Theta(\mu-\wght(I_l))})=1.       
   \end{align}
   for all $l=1,\dots,k$.
 \end{lemma}
 {\bf Proof:}
   Recall \cite[Proposition 4.2.8 (iii)]{b-Joseph} that the right module $V(\mu)^\ast$
   is generated by the highest weight vector $v_\mu^\ast$ and relations 
   \begin{align*}
     v_\mu^\ast y_i=0, \qquad v_\mu^\ast x_i^{2(\mu_i,\alpha_i)/(\alpha_i,\alpha_i)+1}=0,
     \qquad
     v_\mu^\ast t_i =q^{(\mu,\alpha_i)}v_\mu^\ast.
   \end{align*}  
   Using this presentation the first statement of the lemma and the left hand side of
   (\ref{dim1ws1}) are proved by induction on $k$.
   The right hand side of (\ref{dim1ws1}) follows from the fact that there is a nondegenerate
   pairing of $\uqgc$-modules between $V(\mu)$ and $V(-w_0\mu)$. 
   
   Recall \cite[Section 7]{MSRI-Letzter} that $\Theta=-w_0'd$ where $d$ is a 
   diagram automorphism. Hence for any $\mu\in \Ppiplus$ one has
   $\Theta(w_0'w_0\mu)\in \Ppiplus$. To verify (\ref{dim1ws2}) note moreover that 
   $w_0'(\mu-\Theta(\mu))=\mu-\Theta(\mu)$ for all $\mu\in \Ppiplus$ and hence
   $\Theta(w_0'w_0\mu)=\mu$ if $\mu\in \Ptheta$.
   Hence it suffices to show that
   \begin{align}\label{V_thetanu}
        \dim(V(\mu)^\ast_{\nu})=\dim(V(\Theta(w_0'w_0\mu))^\ast_{\Theta(\nu)})
   \end{align}  
   for all $\mu\in \Ppiplus$ and $\nu\in P(\pi)$.
   Moreover, it suffices to verify (\ref{V_thetanu}) in the classical situation. 
   By slight abuse of notation we denote the irreducible left $\gfrak$-module of highest
   weight $\mu$ also by $V(\mu)$.
   To verify (\ref{V_thetanu}) note that
   the vectorspace $V(\mu)$ with the $\gfrak$-module structure defined by 
   $g\cdot v=\theta(g)v$ is isomorphic to $V(\Theta(w_0'w_0\mu))$. Indeed,
   the weight vector $v_{w_0'w_0\mu}\in V(\mu)$ of weight $w_0'w_0\mu$ for the 
   old $\gfrak$-module structure is the highest weight 
   vector for the new $\gfrak$-module structure $\cdot$ and it is of weight 
   $\Theta(w_0'w_0\mu)$.
   $\blacksquare$
 
 \subsection{A weight estimate}
 Before making use of the above Lemma we take another close look at the structure of the 
 central element $d_\mu$. Recall (\ref{dmu=vmu+}) and write $v_\mu\in (\adr \uqgc)\tau(2\mu)$ 
 in the form
 \begin{align}\label{dmuformula}
   v_\mu=\sum_{I,J\in \cIm}a_{IJ} (\adr x_I)(\adr y_J)\tau(2\mu).
 \end{align}
 \begin{lemma}\label{aIJvan}
   Assume $\pt=\emptyset$ and $\mu\in \Ptheta$. Then the following relation holds
   \begin{align}\label{le2mutil}
     v_\mu\in \sum_{-\Theta(\zeta)+\xi\le 2\mutil} 
     (\adr \uqnp_\zeta)(\adr \uqnm_{-\xi})\tau(2\mu).
   \end{align}
   In particular, the coefficients $a_{IJ}$ in (\ref{dmuformula}) vanish if 
   $\wght(J)-\Theta(\wght (I))\not\le 2\mutil$.
 \end{lemma}
 {\bf Proof:} The relation $d_\mu\in \sum_{\nu\le 2\mutil} \uqgc_{-\nu}$ together with
 $d_\mu\in \Bc=\sum B_J T'_\Theta$ implies $d_\mu\in \sum_{\wght(J)\le 2\mutil}B_J T'_\Theta$. 
 Hence, using the explicit form (\ref{Bidef}) of the generators $B_i$, one obtains 
 $d_\mu \in \sum_{-\Theta(\zeta)+\xi\le 2\mutil} \uqnp_\zeta \Tc \uqnm_{-\xi}$ which by
 \cite[Lemma 2.2]{a-KLp} and (\ref{dmu=vmu+}) implies (\ref{le2mutil}).
 $\blacksquare$
 
 \subsection{Yet another projection}
  In the proof of the following lemma and subsequent arguments we will make use of a projection 
 \begin{align*}
   \pi_{IJ}:(\adr \uqgc)\tau(2\mu)\rightarrow  \Cc (\adr x_{I})(\adr y_{J})\tau(2\mu)
 \end{align*}
 onto the space spanned by $(\adr x_{I})(\adr y_{J})\tau(2\mu)$.
 More precisely, for any element 
 $a=\sum_{I',J'\in \cIm}a_{I'J'}(\adr x_{I'})(\adr y_{J'})\tau(2\mu)$ of 
 $(\adr \uqgc)\tau(2\mu)$ and any $I,J\in \cIm$ we define
 $\pi_{IJ}(a)=a_{IJ}(\adr x_{I})(\adr y_{J})\tau(2\mu)$. 
  Note that the formulas (\ref{adrxC}) and (\ref{adryC}) imply that
 \begin{align}\label{weightbound}
   \kow((\adr x_{I})(\adr y_{J})\tau(2\mu))\subseteq 
     \sum_{\makebox[0cm]{$\wght(I')\le \wght(I) \atop \wght(J')\le \wght(J)$}}
     \Cc(\adr x_{I'})(\adr y_{J'})\tau(2\mu))\ot \uqgc.
 \end{align}
 
 We now show that there are sufficiently many nonvanishing coefficients for which the 
 bound of Lemma \ref{aIJvan} is attained.
 We use the notation of Lemma \ref{IinI} and assume that $\mu\in \Ptheta$.
 By (\ref{dim1ws2}) for any $I_l$, $l{=}1,\dots,k$, there exists a uniquely 
 determined $\Ibar_l\in \cIm$ such that $\mu-\wght(\Ibar_l)=\Theta(\mu-\wght(I_l))$.
 Note that by definition $\wght(I_l)-\Theta(\wght(\Ibar_l))=2\mutil$.

 \begin{lemma}\label{nonzero}
   Assume $\pt=\emptyset$ and $\mu\in \Ptheta$ and let $I=(i_1,\dots,i_k)$ be a multiindex 
   satisfying 1) and 2) of Lemma \ref{IinI}.
   For any $l{=}1,\dots,k$ the coefficient $a_{\Ibar_lI_l}$ in (\ref{dmuformula}) does not 
   vanish.
 \end{lemma}
 {\bf Proof:} 
 We proceed by induction on $l$. Note that for the empty multiindex $I_0=()$ relation
 (\ref{dmuX}) implies $a_{\Ibar_0I_0}\neq 0$. 
 Assume now that $a_{\Ibar_{l-1}I_{l-1}}\neq 0$ for some $l$. 
 Recall that by definition $B_i=y_it_i+t^{-1}_{p(i)}x_{p(i)}t_i +s_it_i$ for any 
 $\alpha_i\in \pi$.
 We will use the relation
 \begin{align}\label{Bi-inv}
   (\adr B_{i_l})v_\mu =s_{i_l} v_\mu
 \end{align}
 to show that $a_{\Ibar_{l}I_{l}}\neq 0$.  
 Recall the definition of the coefficients $a_{IJ}$ in (\ref{dmuformula}) and note that
 \begin{align}
   (\adr y_{i_l}t_{i_l})(v_\mu)&\in (\adr t_{i_l})\sum_{I,J\in \cIm}a_{IJ}
   (\adr x_I)(\adr y_{i_l})(\adr y_J)\tau(2\mu) + 
   \bigcap _{\makebox[1cm]{\small$\wght(J)-\Theta(\wght(I))\atop =2\mutil+\alpha_{i_l}$}} 
   \ker(\pi_{IJ})\nonumber\\
   &\subseteq 
   \sum_{\makebox[2cm]{$\wght(J)-\Theta(\wght(I))\atop =2\mutil+\alpha_{i_l}$}} 
   a'_{IJ}  (\adr x_I)(\adr y_J)\tau(2\mu) + 
   \bigcap _{\makebox[1cm]{\small$\wght(J)-\Theta(\wght(I))\atop =2\mutil+\alpha_{i_l}$}} 
   \ker(\pi_{IJ})\label{adyest}
 \end{align}
 for some $a_{IJ}'\in \Cc$. Note moreover that in view of conditions 1) and 2) of Lemma
 \ref{IinI} the assumption $a_{\Ibar_{l-1}I_{l-1}}\neq 0$ 
 implies $a'_{\Ibar_{l-1}I_{l}}\neq 0$. 
 Similarly one calculates
 \begin{align}
   (\adr t^{-1}_{p(i)}x_{p(i)}t_i)(v_\mu)&\sum_{I,J\in \cIm}a_{IJ}
   (\adr t^{-1}_{p(i)}x_{p(i)}t_i) (\adr x_I)(\adr y_J)\tau(2\mu)\nonumber\\
   &= \sum_{I,J\in \cIm}a''_{IJ} (\adr x_I)(\adr y_J)\tau(2\mu)\label{adxest}
 \end{align}
 for some $a''_{IJ}\in \Cc$. Note here, that $a''_{\Ibar_{l-1}I_{l}}\neq 0$ if and only if
 $a_{\Ibar_{l}I_{l}}\neq 0$. 
 Finally by Lemma \ref{aIJvan} both $v_\mu$ and $(\adr t_i)(v_\mu)$ are contained in 
 \begin{align*}
   \bigcap _{\makebox[1cm]{\small$\wght(J)-\Theta(\wght(I))\atop =2\mutil+\alpha_{i_l}$}} 
   \ker(\pi_{IJ})
 \end{align*}
 In view of (\ref{Bi-inv}) the relations (\ref{adyest}) and (\ref{adxest}) imply that 
 $a'_{\Ibar_{l-1}I_{l}}=-a''_{\Ibar_{l-1}I_{l}}$. Now 
 $a_{\Ibar_{l}I_{l}}\neq 0$ follows from the properties of $a'_{\Ibar_{l-1}I_{l}}$ and
 $a''_{\Ibar_{l-1}I_{l}}$ noted above. 
 $\blacksquare$
 
 \subsection{Finding $B_i\in L_{-2w_0\mu}\Ttcp$}\label{finalstep}
 We are now in a position to prove that for any $\mu\in \Ptheta$ the right coideal
 $L_{-2w_0\mu}$ contains sufficiently many generators of $\Bc$.
 \begin{proposition}\label{finalprop}
   Assume that $\pt=\emptyset$ and $\mu\in \Ptheta$. Then for any 
   $\alpha_i\in \pi$ there exists $t(i)\in \Ttcp$ such that $B_it(i) \in L_{-2w_0\mu}$.
 \end{proposition}
 \noindent{\bf Proof:}
 As $\pt=\emptyset$ we can choose $I=(i_1,\dots,i_{k})$ with $i_k=i$ such that
 assumptions 1) and 2) of Lemma \ref{IinI} hold. Note that 
 $\wght(I_{k-1}){-}\Theta(\wght(\Ibar_k))=2\mutil{-}\alpha_i$.
 In view of (\ref{dmuformula}) and Lemma \ref{aIJvan} the relation (\ref{weightbound}) 
 implies
 \begin{align}
   (\pi_{\Ibar_k I_{k-1}}\ot \id)\kow (v_\mu)\in (\pi&_{\Ibar_k I_{k-1}}\ot \id)\kow 
   (\pi_{\Ibar_k I_k}(v_\mu))\nonumber\\
   & + (\adr x_{\Ibar_k})(\adr y_{I_{k-1}})\tau(2\mu)\ot 
   \sum_{\nu>-\alpha_i} \uqgc_{\nu}. \label{uptolowerweights}
 \end{align}
 Using formulas (\ref{adrxC}) and (\ref{adryC}) and assumptions 1) and 2) of Lemma \ref{IinI} again one obtains
 \begin{align}\label{lowestpart}
   (\pi&_{\Ibar_k I_{k-1}}\ot \id)\kow (\pi_{\Ibar_k I_k}(v_\mu))\nonumber\\ &=
   a_{\Ibar_k I_k}(\adr x_{\Ibar_k})(\adr y_{I_{k-1}})\tau(2\mu)\ot \tau(-\wght(\Ibar_k))
   (t y_i {-} q^{-a}y_i t)
 \end{align}
 where $a:=(\alpha_i,\alpha_{i_1}+\dots+\alpha_{i_{k-1}})$ and 
 $t:=\tau(2\mu-\alpha_{i_1}-\dots-\alpha_{i_{k-1}})$.
 Note that if $(\mu,\alpha_i)\neq 0$ then by construction $k=1$ and hence $a=0$ and 
 $t=\tau(2\mu)$. If on the other hand $(\mu,\alpha_i)= 0$ then 
 $a=(\alpha_i,\alpha_{i_{k-1}})\neq 0$. One obtains
 \begin{align}\label{ty-yt}
   ty_i-q^{-a}y_it=
     \begin{cases} (q^{-(2\mu,\alpha_i)}{-}1) y_i\tau(2\mu) &\mbox{ if $(\mu,\alpha_i)\neq 0$,}\\
                   (q^a{-}q^{-a}) y_i t &\mbox{ if $(\mu,\alpha_i)=0$.}      
     \end{cases}
 \end{align}
 In any case, (\ref{uptolowerweights}), (\ref{lowestpart}), (\ref{ty-yt}), and Lemma 
 \ref{nonzero} imply that there exists
 \begin{align*}
   b\in L_{-w_0\mu}\cap (y_i t \tau(-\wght(\Ibar_k)) +\sum_{\nu\ge 0}\uqgc_{\nu})
 \end{align*}
 Hence $B_i t(i)+m\in L_{-w_0\mu}$ for some $m\in \Cc \Ttcp$ and
 $t(i)=\tau(2\mu{-}\wght(I_k){-}\wght(\Ibar_k))$. 
 Note that the relation $\wght(I_k)-\Theta(\wght(\Ibar_k))=2\mutil$ implies
 $t(i)=\tau(\mu-\wght(I_k)+\Theta(\mu-\wght(I_k)))\in \Ttcp$.
 Applying the coideal 
 property of $L_{-w_0\mu}$ to the element $B_i t(i)+m$ one obtains $B_it(i)\in L_{-w_0\mu}$.
 $\blacksquare$

\begin{corollary}
  If $\pt=\emptyset$ then Conjecture \ref{conj2} holds.
\end{corollary}
\noindent {\bf Proof:} By the above Proposition the coideal $L_\mu$ contains elements 
    $B_it(i)$ for all $\alpha_i\in \pi$ and some $t(i)\in \Ttcp$. If $\vep(B_i)=0$ then any
    element in $V(\nu)^{L_\mu}$ is also invariant under $B_i$.
    Note that $\vep(B_i)= 0$ holds for all $\alpha_i$ such that 
    $\Theta(\alpha_i)=-\alpha_{p(i)}\neq -\alpha_i$.
    By \cite[Theorem 7.1(iii)]{a-Letzter03} one has
    \begin{align*}
      B_iB_{p(i)}-B_{p(i)}B_i=
      \frac{t_it_{p(i)}^{-1}-t_i^{-1}t_{p(i)}}{q^{(\alpha_i,\alpha_i)/2}-q^{-(\alpha_i,\alpha_i)/2}}
    \end{align*}
    and hence $B_i$, $B_{p(i)}$, $t_it_{p(i)}^{-1}$, and $t_i^{-1}t_{p(i)}$ generate a 
    subalgebra isomorphic to $U_q(\mathfrak{sl}_2)$. Therefore invariance under $B_i$ and 
    $B_{p(i)}$ implies invariance under $t_it_{p(i)}^{-1}$.
    Hence $V(\nu)^{L_\mu}\subseteq V(\nu)^{\Ttcp}$. Proposition \ref{finalprop} now implies
    that any element in $V(\nu)^{L_\mu}$ is also invariant under all $B_i$.
    $\blacksquare$

\providecommand{\bysame}{\leavevmode\hbox to3em{\hrulefill}\thinspace}
\providecommand{\MR}{\relax\ifhmode\unskip\space\fi MR }
\providecommand{\MRhref}[2]{%
  \href{http://www.ams.org/mathscinet-getitem?mr=#1}{#2}
}
\providecommand{\href}[2]{#2}

\textsc{Stefan Kolb, Mathematics Department, Virginia Polytechnic Institute
         and State University, Blacksburg, VA 24061, USA.}
         
 \textit{E-mail address:} \texttt{kolb@math.vt.edu}

\end{document}